%% file: Revision.tex
\documentclass[11pt,reqno]{amsart}
\usepackage[margin=1.1in]{geometry}

\usepackage{amsmath,amsthm,amssymb}
\usepackage{mathrsfs} 
\usepackage{mathabx}
\usepackage{array} 
 
\usepackage{enumitem}

\usepackage{tikz}
\usepackage{pgfplots}

\usetikzlibrary{arrows}
\usepackage[unicode]{hyperref}

\usepackage{bm} 
\usepackage{amsmath, amsfonts, amssymb, amsthm} 
\usepackage{mathrsfs}
\usepackage{graphicx} 
\usepackage{float} 
\usepackage{enumitem} 
\usepackage{fancyhdr} 
\usepackage{caption}  
\usepackage{booktabs} 
\usepackage[comma,sort&compress, numbers]{natbib}
\usepackage{subcaption}

\numberwithin{equation}{section} 

\newtheorem{theorem}{Theorem}[section]
\newtheorem{cor}[theorem]{Corollary}

\newtheorem{lemma}[theorem]{Lemma} 
\theoremstyle{definition} 
\newtheorem{definition}[theorem]{Definition}
\newtheorem{remark}[theorem]{Remark}
\newtheorem{example}[theorem]{Example}

\def\R{\mathbb{R}}

\def\P{\mathbb{P}}
\def\E{\mathbb{E}}

\def\Var{\mathrm{Var}}
\def\Cov{\mathrm{Cov}}

\def\1{\mathbbm{1}}

\newcommand{\cQ}{\bm{Q}}
\newcommand{\cR}{\bm{R}}

\DeclareMathOperator{\ext}{\mathcal{E}}

\allowdisplaybreaks

\title[Thresholds and Fluctuations of Submultiplexes in Random Multiplex Networks]{Thresholds and Fluctuations of Submultiplexes in Random Multiplex Networks }

\author[Bhattacharya]{Bhaswar B. Bhattacharya} 
\address{Department of Statistics and Data Science, University of Pennsylvania, Philadelphia, USA} \email{bhaswar@wharton.upenn.edu}

\author[Bhowal]{Sanchayan Bhowal}
\address{Department of Statistics, Stanford University, California, USA} 
\email{sbhowal@stanford.edu}

\author[Das]{Karambir Das}
\address{Indian Institute of Science, Bengaluru, India} 
\email{karambirdas@iisc.ac.in} 

\author[Eslava]{Laura Eslava}
\address{Probability and Statistics Department, IIMAS-UNAM, Mexico City, Mexico} 
\email{laura@sigma.iimas.unam.mx} 

\author[Karmakar]{Shaibal Karmakar}
\address{International Centre for Theoretical Sciences, Bengaluru, India} 
\email{shaibal.karmakar@icts.res.in} 

\begin{document}
\begin{abstract}
    In a multiplex network, a common set of nodes is connected through different types of interactions, each represented as a separate graph (layer) within the network. In this paper, we study the asymptotic properties of submultiplexes, the counterparts of subgraphs (motifs) in single-layer networks, in the correlated Erd\H{o}s-R\'{e}nyi multiplex model. This is a random multiplex model with two layers, where the graphs in each layer marginally follow the classical (single-layer) Erd\H{o}s-R\'{e}nyi model, while the edges across layers are correlated.
    We derive the precise threshold condition for the emergence of a fixed submultiplex $\bm{H}$ in a random multiplex sampled from the correlated Erd\H{o}s-R\'{e}nyi model. Specifically, we show that the satisfiability region, the regime where the random multiplex contains a diverging number of copies of $\bm{H}$, is a polyhedral region in $\mathbb{R}^3$. Furthermore, within this region the count of $\bm{H}$ is asymptotically normal, with an explicit convergence  rate in the Wasserstein distance. We also establish various Poisson approximation results for the count of $\bm{H}$ on the boundary of the threshold, which depends on a notion of {\it balance} of submultiplexes. Collectively, these results provide an asymptotic theory for small submultiplexes in the correlated multiplex model, analogous to the classical theory of small subgraphs in random graphs.
\end{abstract}
\maketitle

\section{Introduction}

The basic model of a random graph is the Erd\H{o}s-R\'enyi model $G(n, p)$, which is a graph with vertex set $[n]:=\{1, 2, \ldots, n\}$, where each edge $(i, j)$ is present independently with probability $p \in [0, 1]$, for $1 \leq i < j \leq n$. This model forms the foundation of random graph theory and network analysis, and its properties have been extensively studied. In this paper, we consider {\it multiplex networks}, which consist of multiple graphs defined on a common set of nodes. In such networks, the edges of each graph (referred to as the layers of the multiplex) represent different types of interactions. For example, in social networks the
different layers can correspond to types of relationships (such as friendship, collaboration, or family ties) or modes of communication (such as email, phone, or message) \cite{wasserman1994social,fienberg1985statistical,pattison1999logit,florentinegraph}; in transportation networks, layers represent
different modes of transportation  \cite{cardillo2013emergence,strano2015multiplex,bergermann2021multiplex,multiplexindex} (also known as multi-mode graphs \cite{kirkpatrick2025shortest}); and in biological networks, such as those describing connections in the nervous system, layers may represent synaptic (chemical) or electrical links between neurons \cite{bentley2016multilayer}; to name a few. The study of multiplex networks has been growing rapidly in recent years, with applications appearing across a wide range of disciplines (see the monographs \cite{bianconi2018multilayer,networksbiology,multiplexnetworkproperties} and the references therein).

A basic model of a random multiplex is the {\it correlated Erd\H{o}s-R\'enyi model}. In this model, each layer marginally follows the classical (single-layer) Erd\H{o}s-R\'enyi model, while the edges across the layers are dependent. We define this model formally below, considering two layers for simplicity. For this, we need to introduce a few notations. A multiplex with two layers will be denoted by $\bm{G} = (V(\bm{G}), G^{(1)}, G^{(2)})$, where $V(\bm{G})$ is the common vertex set and  $G^{(1)} = (V(\bm G), E(G^{(1)}))$ and $G^{(2)} = (V(\bm G), E(G^{(2)}))$ are the two layers, which are two simple graphs with vertex set $V(\bm G)$ and edge sets $E(G^{(1)})$ and $E(G^{(2)})$, respectively. The union of the edge sets across the two layers will be denoted by $E(\bm{G})= E(G^{(1)}) \cup E(G^{(2)})$.

\begin{definition}\label{definition:correlated} In the {\it correlated Erd\H{o}s-R\'enyi} model a random multiplex $\bm G_n := ([n], G_n^{(1)}, G_n^{(2)})$ is generated as follows: Let $(a_{ij}^{(1)}, a_{ij}^{(2)})_{1 \leq i < j \leq n}$ be an i.i.d. collection of pairs of correlated Bernoulli random variables distributed as: 
$$\P(a_{ij}^{(1)}=1) = p_1, ~ \P(a_{ij}^{(2)}=1) = p_2,  \text{ and } \P(a_{ij}^{(1)}=1, a_{ij}^{(2)}=1) = p_{12} ,$$
    where
    \begin{align}
        \label{eq:p}
        (p_1, p_2, p_{12}) \in (0, 1)^3 \text{ is such that } \max\{ 0, p_1 + p_2 - 1 \} \leq p_{12} \leq \min\{p_1, p_2\}.
    \end{align}
     Then, for $1 \leq i < j \leq n$, the edge $(i,j)$ is included in the two layers as follows: 
   \begin{itemize}
\item if $a_{ij}^{(1)} = 1$ and $a_{ij}^{(2)} = 0$, then $(i, j) \in E(G_n^{(1)})$ and $(i, j) \notin  E(G_n^{(2)})$, 
\item if $a_{ij}^{(1)} = 0$ and $a_{ij}^{(2)} = 1$, then $(i, j) \notin E(G_n^{(1)})$ and $(i, j) \in  E(G_n^{(2)})$, 
\item if $a_{ij}^{(1)} = a_{ij}^{(2)} = 1$, then $(i, j) \in E(G_n^{(1)}) \cap  E(G_n^{(2)})$,   
\item if $a_{ij}^{(1)} = a_{ij}^{(2)} = 0$, then $(i, j) \notin E(G_n^{(1)}) \cup  E(G_n^{(2)})$. 
\end{itemize} 
    (Note that the graphs $G_n^{(1)}$ and $G_n^{(2)}$ are marginally distributed as Erd\H{o}s-R\'enyi models $G(n, p_1)$ and $G(n, p_2)$, respectively, and the edges appear in both layers with probability $p_{12}$.) We write $\bm G_n\sim G(n,p_1,p_2,p_{12})$ for a random multiplex generated as above, and refer to
$G(n,p_1,p_2,p_{12})$ as the correlated Erd\H{o}s--R\'enyi multiplex model. Figure \ref{fig:randommultiplex} (a) shows a  sample from the correlated Erd\H{o}s-R\'enyi model with $n=15$ and $p_1=p_2=0.15$, $p_{12}=0.05$.    
\end{definition}

\begin{figure}[h]
    \begin{subfigure}[c]{0.5\linewidth}
        \begin{center}
            \vspace{-0.15in}
            \input{plots/GraphMultiplex.tex} \\
            \vspace{0.1in}
            { (a) }
        \end{center}
    \end{subfigure}
    \hspace{2cm}
    \begin{subfigure}[c]{0.18\linewidth}

        \begin{center}
            \begin{tikzpicture}[thick, scale=0.725, every node/.style={draw, circle, fill=gray!30, inner sep=1pt, minimum size=8mm}]

                \node (1) at (0, 0) {1};
                \node (2) at (3, 0) {2};
                \node (3) at (1.5, 2.5) {3};

                \draw[blue, thick] (1) -- (3);
                \draw[blue, thick] (3) -- (2);
                \draw[red, thick] (2) -- (1);

                \draw[red,  thick] (2) to[bend right=30] (3);

                \node[draw=none, fill=none, text width=1cm] at (1.05, 3.75) {$\bm{H}$};

            \end{tikzpicture} \\
            { (b) }
        \end{center}
    \end{subfigure}
    \caption{ \small{ (a) A sample $\bm{G}_n$ from the correlated Erd\H{o}s-R\'enyi multiplex with $n=15$ and $p_1=p_2=0.15$, $p_{12}=0.05$. The edges in layer $1$ are colored blue, while those in layer $2$ are colored red. (b) A fixed multiplex $\bm{H}$ on three vertices. $X(\bm H, \bm G_n)$ counts the number of copies of $\bm H$ in $\bm G_n$ (see \eqref{eq:NHGnM} for the formal definition). } }
    \label{fig:randommultiplex}
\end{figure}

The correlated Erd\H{o}s-R\'enyi model initially appeared in a different guise in the context of network de-anonymization \cite{pedarsani2011privacy}. Recently, this has been adopted as the canonical model in random graph matching and related problems (see \cite{mao2024testing,mao2023random,lyzinski2015graph,lyzinski2014seeded,huang2025correlation,barak2019efficient}, among several others). In this context, a prototypical statistical question is to detect whether two networks are edge-correlated through some latent vertex correspondence. Test statistics for such problems often involve the covariance between subgraph count statistics across the two networks \cite{mao2024testing,huang2025correlation,barak2019efficient}. Analyzing these quantities requires understanding the asymptotic properties of counts of two overlapping subgraphs defined on a shared set of nodes, which naturally connects to the study of submultiplexes within a multiplex network.

Motivated by the above problems and the growing need for the multiplex framework in modeling large-scale complex systems, in this paper we aim to understand the asymptotic behavior of submultiplexes. Analogous to the role of subgraphs (motifs) in single-layer networks, submultiplexes are key local features that capture important structural properties in multiplex networks. Focusing on the correlated Erd\H{o}s-R\'enyi model, we begin by deriving the threshold for the emergence (containment) of a fixed submultiplex $\bm H$ in a realization of the random multiplex $\bm G_n \sim G(n, p_1, p_2, p_{12})$. Specifically, we determine the precise threshold condition under which $X(\bm H, \bm G_n)$, the number of copies of $\bm H$ in $\bm G_n$ (see Section \ref{sec:multiplex} for the formal definition and Figure \ref{fig:randommultiplex} for an illustration), converges in probability to zero or to infinity. The threshold condition corresponds to a minimization problem that is determined by the submultiplexes of $\bm{H}$ with the smallest expected number of copies. By parametrizing the model so that the connection probabilities $p_1, p_2, p_{12}$ decay polynomially, we show that the satisfiability region, that is, the regime where $\bm{G}_n$ contains many copies of $\bm H$, is a polyhedral subset of $\mathbb R^3$. Next, we show that whenever $\bm G_n$ contains many copies of $\bm H$, that is, for all parameter values in the satisfiability region, $X(\bm H, \bm G_n)$ is asymptotically normal after appropriate normalization. Additionally, our result provides an explicit convergence rate for this normal approximation in the Wasserstein distance (see Theorem \ref{thm:ZHGn}).
We next analyze the regime in which the parameter values lie on the boundary of the satisfiability region (the threshold surface). At such boundary points, the asymptotic behavior of $X(\bm  H, \bm  G_n)$ depends on the structure of the minimizers of the threshold condition. If the set of minimizers contains $\bm H$ itself (in which case we say $\bm  H$ is {\it balanced} for that parameter value), then the distribution of $X(\bm  H, \bm  G_n)$ depends on whether this minimizer is unique. In the {\it strictly balanced} case, where the minimizer is achieved uniquely at $\bm  H$, we show that $X(\bm  H, \bm  G_n)$ converges in distribution to a Poisson, under an additional nondegeneracy condition (see Theorem~\ref{thm:balancedH}). On the other hand, if all minimizers occur at proper submultiplexes of $\bm  H$, then an appropriate rescaling reduces the problem to the balanced case (see Theorem~\ref{thm:unbalancedH}). Together, these results develop an analogous theory for small submultiplexes in the correlated multiplex model, paralleling the classical results on small subgraphs in Erd\H{o}s-R\'{e}nyi random graphs. Finally, in Section \ref{sec:summary} we list a few questions and directions for future research.

\subsection*{Asymptotic Notation}

Throughout the paper we will use the following asymptotic notations: For two sequences $a_n$ and $b_n$ we will write $a_n \lesssim b_n$ if for all $n$ large enough $a_n \leq C _1 b_n$, for some constant $C_1 > 0$. Similarly, $a_n \gtrsim b_n$ will mean  $a_n \geq C_2 b_n$ and $a_n \asymp b_n$ will mean $C_2 b_n \leq a_n \leq C_1 b_n$, for $n$ large enough and constants $C_1, C_2 > 0$. Subscripts in the above notation, for example $\lesssim_{\square}$ and $\asymp_{\square}$ denote that the hidden constants may depend on the subscripted parameters. Moreover, $a_n \ll b_n$, $a_n\gg b_n$, and $a_n \sim b_n$ mean $a_n/b_n \rightarrow 0$, $a_n/b_n \rightarrow \infty$, $a_n/b_n \rightarrow 1$, as $n \rightarrow \infty$.

\section{Statements of Results}
\label{sec:results}


\subsection{Thresholds for Submultiplex Containment}
\label{sec:multiplex}

Thresholds for the emergence of fixed subgraphs in the random
graph $G(n, p)$ were among the first problems studied by Erd\H{o}s and R\'enyi in their fundamental paper \cite{erdos1960evolution}. Bollob\'as \cite{bollobas1981threshold} resolved this problem in its full generality a few years later, and this is now a textbook result in random graph theory (see, for example \cite[Chapter 3]{janson2011random}). Our first aim in this paper is to derive the analogous result for submultiplexes in the  correlated Erd\H{o}s-R\'enyi multiplex model. We begin by formalizing the notion of containment in a multiplex network. For two multiplexes $\bm F= (V(\bm F), F^{(1)}, F^{(2)}) $ and $\bm G = (V(\bm G), G^{(1)}, G^{(2)})  $, we write $\bm F\subseteq\bm G$ and say that $\bm F$ is a \emph{submultiplex} of $\bm G$ if $V(\bm F)\subseteq V(\bm G)$
and $E(F^{(i)})\subseteq E(G^{(i)})$, for $i \in \{1, 2\}$. Given a fixed multiplex $\bm H = (V(\bm{H}), H^{(1)}, H^{(2)})$  and an injective map $\phi: V(\bm H) \rightarrow [n]$, denote by $\bm H_{\phi}$ the image multiplex with vertex set $\phi(V(\bm{H}))$ and edge sets $\phi(E(H^{(i)})) = \{(\phi(u), \phi(v)): (u, v) \in E(H^{(i)}) \}$, for $i \in \{1, 2\}$. Then the number of (unlabeled) copies of $\bm H$ in $\bm G_n = ([n], G_n^{(1)}, G_n^{(2)}) \sim G(n, p_1, p_2, p_{12})$, generated from the correlated Erd\H{o}s-R\'enyi model as in Definition \ref{definition:correlated}, is defined as: 
\begin{equation}\label{eq:NHGnM}
X(\bm H,\bm G_n) = \sum_{\overline{\bm H}'\in\mathcal C_n(\bm H)}
\bm 1 \{\overline{\bm H}'\subseteq\bm G_n \} , 
\end{equation}
where $\mathcal C_n(\bm H) := \left\{ \bm H_{\phi}: \phi: V(\bm H) \rightarrow [n] \text{ is injective} \right\}$. The following result gives the threshold for the emergence of $\bm{H}$ in $G(n, p_1, p_2, p_{12})$. Throughout we are going to assume $p_1, p_2$ (and as a result $p_{12}$) are bounded away from 1.


\begin{theorem}
    \label{thm:threshold}
    Suppose $\bm H = (V(\bm{H}), H^{(1)}, H^{(2)})$ is a fixed multiplex with at least one edge (that is, $|E(\bm{H})| \geq 1$). Then for $X(\bm{H}, \bm G_n)$ as defined above, the following holds:
    \begin{equation}\label{eq:XHGn}
        X(\bm{H}, \bm G_n) \stackrel{P} \rightarrow \begin{cases}
            0      & \Phi_{\bm{H}} \ll 1  , \\
            \infty & \Phi_{\bm{H}} \gg 1  ,  
        \end{cases}  
    \end{equation}
    where
    \begin{align}\label{eq:PhiH}
        \Phi_{\bm{H}} & =   \Phi_{\bm{H}} (n, p_1, p_2, p_{12}) \nonumber                         \\
                      & :=\min_{\substack{\bm{F} = (V(\bm{F}), F^{(1)}, F^{(2)}) \subseteq \bm{H} \\ |E(\bm{F})| \geq 1 }} n^{|V(\bm{F})|}p_1^{|E(F^{(1)})\setminus E(F^{(2)})|} p_2^{|E(F^{(2)})\setminus E(F^{(1)}) |} p_{12}^{|E(F^{(1)})\cap E(F^{(2)}) |} .
    \end{align}
\end{theorem}

The proof of Theorem \ref{thm:threshold} relies on standard applications of the first and second moment methods (see Section \ref{sec:thresholdpf}).
From the proof it will be evident that
$$\Phi_{\bm{H}} \asymp_{\bm H} \min\{ \E [X(\bm{F}, \bm G_n) ] : \bm{F} \subseteq \bm{H},  |E(\bm{F})| \geq 1 \}.$$
In other words, the threshold for the emergence of a fixed submultiplex is determined
by the rarest (in terms of expected number of copies) submultiplex of $\bm{H}$. The same phenomenon determines the threshold for the appearance of any fixed subgraph in an Erd\H{o}s-R\'enyi random graph (see \cite[Theorem 3.4]{janson2011random}).

\begin{remark}\label{remark:graphdistribution}
    One can recover the threshold for
    subgraphs in Erd\H{o}s-R\'{e}nyi random graphs by considering one of the layers in $\bm{H}$ to be empty. Specifically, if $H^{(2)}$ is the empty graph and $H^{(1)}= H = (V(H), E(H))$ is a fixed simple graph, then \eqref{eq:PhiH} simplifies to
    $$\Phi_{H} =\min_{\substack{F \subseteq H, |E(F)| \geq 1 }} n^{|V(F)|}p_1^{|E(F)|}
        .$$
    In this case, the condition $\Phi_{H}  \gg 1$ is equivalent to $ n p_1^{m_H} \gg 1 $, where
    \begin{align}\label{eq:mH}
        m_H:=\max_{F\subseteq H,\,|V(F)|>0}\frac{|E(F)|}{|V(F)|} .
    \end{align}
    This matches the threshold for the emergence of $H$ in the Erd\H{o}s-R\'enyi graph $G(n, p_1)$ (see \cite[Theorem 3.4]{janson2011random}).    If the maximum in \eqref{eq:mH} is attained by $H$ itself, that is, $m(H)=\frac{|E(H)|}{|V(H)|}$, then $H$ is said to be a {\it balanced} graph. Otherwise, $H$ is said to be an {\it unbalanced} graph. Moreover, if the maximum in \eqref{eq:mH} is uniquely attained at $H$, then $H$ is called a {\it strictly balanced} graph.
\end{remark}

To better visualize the phase transition in Theorem \ref{thm:threshold} it is instructive to parameterize the probabilities as:
\begin{align}\label{eq:exponent}
    p_1 = n^{-\theta_1}, \quad p_2 = n^{-\theta_2}, \quad \text{ and } \quad p_{12} = n^{-\theta_{12}}.  
\end{align}
The conditions in \eqref{eq:p} impose the following restrictions on the values of $\theta_1, \theta_2, \theta_{12}$: 
$$\Theta := \{ \bm \theta = (\theta_1, \theta_2, \theta_{12} ) : \theta_1, \theta_2, \theta_{12} > 0 \text{ and } \theta_{12} \geq \max\{\theta_1, \theta_2\}\}.$$
Then \eqref{eq:XHGn} implies: 
\begin{equation}\label{eq:polyregion}
    X(\bm{H}, \bm G_n) \stackrel{P} \rightarrow \begin{cases}
        0      & \Delta_{\bm H} < 0 \\
        \infty & \Delta_{\bm H} > 0
    \end{cases} ,
\end{equation}
where
\begin{align}\label{eq:parameterH}
    \Delta_{\bm H} & = \Delta_{\bm{H}}(\bm \theta)     \\
                   & = \min_{ \substack{\bm{F} \subseteq \bm{H} \\ |E(\bm{F})| \geq 1 } } \left\{ |V(\bm{F})|- \theta_1 |E(F^{(1)})\setminus E(F^{(2)})|- \theta_2 |E(F^{(2)})\setminus E(F^{(1)})|- \theta_{12} |E(F^{(1)})\cap E(F^{(2)})| \right \} ,  \nonumber
\end{align}
with the minimum taken over all submultiplexes $\bm{F} = (V(\bm{F}), F^{(1)}, F^{(2)}) \subseteq \bm{H}$.
Note that condition $\Delta_{\bm H} > 0$ is equivalent to the condition for all $\bm{F} \subseteq \bm{H}$ with $|E(\bm F)| \geq 1$,
\begin{align}\label{eq:VF}
    \ell_{\bm F}(\bm \theta) & := |V(\bm{F})|- \theta_1 |E(F^{(1)})\setminus E(F^{(2)})|- \theta_2 |E(F^{(2)})\setminus E(F^{(1)})|- \theta_{12} |E(F^{(1)})\cap E(F^{(2)})|  \nonumber \\
                             & > 0 ,
\end{align}
This defines an open halfspace in $\mathbb R^3$ that contains the origin. Hence,
the parameter values for which the count of $\bm{H}$ in $\bm G_n$ diverges, are given by   (recall \eqref{eq:polyregion}):
$$\partial_{\bm H}^+ := \{(\theta_1, \theta_2, \theta_{12}) \in \Theta: \Delta_{\bm H} > 0\} $$
is formed by the intersection of finitely many halfspaces. Hence, it is a (possibly unbounded) convex polyhedral region contained in the open positive orthant of $\R^3$.
We will refer to this region as the {\it satisfiability region} for $\bm{H}$.
Note that the satisfiability region $\partial_{\bm H}^+ $ is determined by two types of linear constraints:
\begin{itemize}

    \item Constraints of the form \eqref{eq:VF} which are {\it open} halfspaces in $\R^3$. These will be referred to as the {\it main constraints}.


    \item The {\it feasibility constraint} $\theta_{12} \geq \max\{\theta_1, \theta_2\}$ which is the intersection of the closed halfspaces $\theta_{12} \geq \theta_1$ and $\theta_{12} \geq \theta_2$. (The $\theta_1, \theta_2, \theta_{12} > 0 $ are also feasibility constraints, which are enforced by considering the positive orthant of $\R^3$.)

\end{itemize}
Hence, a boundary piece of $\partial_{\bm H}^+$ may or may not be included in the
region, depending on whether it arises from a strict main constraint or a non-strict feasibility constraint. The {\it unsatisfiability region}
$$ \partial_{\bm H}^- := \{(\theta_1, \theta_2, \theta_{12}) \in \Theta: \Delta_{\bm H} <  0\} $$
consists of the parameter values for which there are asymptotically no copies of $\bm{H}$ in $\bm G_n$. These 2 regions are separated by the {\it threshold surface}
$$\partial_{\bm H} := \{(\theta_1, \theta_2, \theta_{12}) \in \Theta: \Delta_{\bm H} =  0\} , $$
which forms a polyhedral surface in $\mathbb R^3$. The submultiplexes of $\bm H$ which determine the condition $\Delta_{\bm H} =  0$ will play an important role in our analysis. This is formalized in the following definition:

\begin{definition}
    Fix a multiplex $\bm H =(V(\bm H), H^{(1)}, H^{(2)})$ and a point $\bm \theta= (\theta_1, \theta_2, \theta_{12}) \in \partial_{\bm H}$. A submultiplex
    $\bm F =(V(\bm F), F^{(1)}, F^{(2)}) \subseteq \bm H$, with $|E(\bm F)| \geq 1$,  is said to be $\bm \theta$-{\it extremal}, if $\ell_{\bm F}(\bm \theta)  = 0$ (recall \eqref{eq:VF}).
    \label{defn:boundaryH}
\end{definition}

To illustrate the ideas described above, we consider the following examples.


\begin{figure*}[h]
    \begin{subfigure}[b]{0.18\linewidth}
        \begin{center}
            \begin{tikzpicture}[thick, scale=0.725, every node/.style={draw, circle, fill=gray!30, inner sep=1pt, minimum size=8mm}]

                \node (1) at (0, 0) {1};
                \node (2) at (3, 0) {2};
                \node (3) at (1.5, 2.5) {3};

                \draw[blue, thick] (1) -- (3);
                \draw[blue, thick] (3) -- (2);
                \draw[blue, thick] (2) -- (1);

                \draw[red,  thick] (2) to[bend right=30] (3);

            \end{tikzpicture} \\
            { (a) }
        \end{center}
    \end{subfigure}
    \hfill
    \begin{subfigure}[b]{0.45\linewidth}
        \begin{center}
            \includegraphics[width=\linewidth]{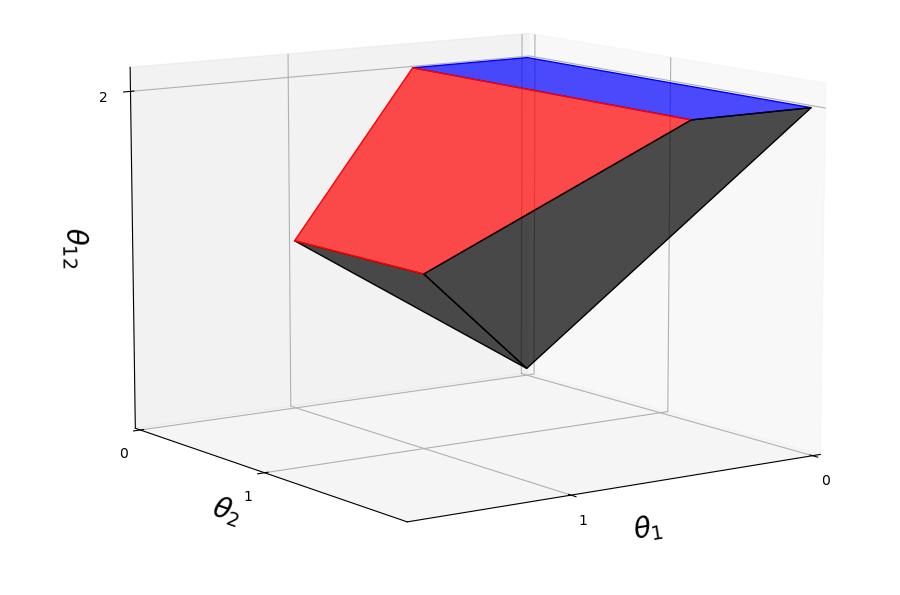} \\
            { (b) }
        \end{center}
    \end{subfigure}
    \hfill
    \begin{subfigure}[b]{0.35\linewidth}
        \begin{center}
            \input{plots/Multiplex_1.tex} \\
            { (c) }
        \end{center}
    \end{subfigure}
    \caption{ (a) The edge-triangle multiplex $\bm{\mathcal{R}}$ defined in \eqref{eq:H123}, (b) the full phase diagram in 3D, and (c) the 2D phase diagram assuming $\theta_1=\theta_2= \theta$. }
    \label{fig:edgetriangle}
\end{figure*}

\begin{example}\label{example:edgetriangle}
    Consider the edge-triangle multiplex as shown in Figure \ref{fig:edgetriangle} (a):
    \begin{align}\label{eq:H123}
        \bm{\mathcal{R}} = (\{1, 2, 3\}, \mathcal{R}^{(1)}, \mathcal{R}^{(2)} ) ,
    \end{align}
    with $E(\mathcal{R}^{(1)})=  \{(1, 2), (2, 3), (1, 3)\}$, which corresponds to the triangle in layer 1 with edges colored in blue, and $E(\mathcal{R}^{(2)}) =  \{(2, 3)\}$, which corresponds to the edge in layer 2 colored in red. In this case, the main constraints determining the boundary of the region $\partial_{\bm{\mathcal{R}}}^{+}$ are:
    \begin{itemize}
        \item $2 \theta_1 + \theta_{12}  <  3 $: This corresponds to taking $\bm F = \bm{\mathcal{R}}$ in \eqref{eq:VF}.
        \item $\theta_{12} < 2$:   This corresponds to taking $\bm F = \bm{\mathcal{O}} := (\{2, 3\}, \mathcal{O}^{(1)}, \mathcal{O}^{(2)})$, where $\mathcal{O}^{(1)}$ and $\mathcal{O}^{(2)}$ are both the single edge $(2,3)$.
    \end{itemize}
    In addition one has the feasibility constraint $\max\{\theta_1, \theta_2\} \leq  \theta_{12}$. These together generate the 3-dimensional satisfiability region shown in Figure \ref{fig:edgetriangle} (b): The red plane corresponds to the equation $2 \theta_1 + \theta_{12}  =3 $, the blue plane corresponds to $\theta_{12} = 2$, and the black planes correspond to $\max\{\theta_1, \theta_2\} =  \theta_{12}$. In Figure \ref{fig:edgetriangle} (c) we show the 2-dimensional slice of this region where $\theta_1=\theta_2= \theta$ (this corresponds to the case where the graphs in the marginal layers have the same distribution). Then $\partial_{\bm{\mathcal{R}}}^+$ is the two-dimensional region:
    $$\partial_{\bm{\mathcal{R}}}^+ := \{ (\theta, \theta_{12}) :  \theta, \theta_{12} > 0 , \theta \leq \theta_{12},  2 \theta + \theta_{12}  <  3, \theta_{12} < 2 \} .$$
    This is the blue polygonal region together with the blue line segment $(O, Q_1)$ (excluding the origin $O$ and the point $Q_1$) shown in Figure \ref{fig:edgetriangle} (c). The threshold curve $\partial_{\bm{\mathcal{R}}}$ is the red polyline segment $[Q_1, Q_2, Q_3)$ (including the points $Q_1$ and $Q_2$ and excluding the point $Q_3$) and the unsatisfiability region $\partial_{\bm{\mathcal{R}}}^{-}$ is shown in grey in Figure \ref{fig:edgetriangle} (c).  Note that for all points $\bm \theta = (\theta, \theta_{12})$ on the open segment $(Q_1, Q_2)$, $\bm{\mathcal{R}}$ is the unique $\bm \theta$-extremal multiplex. On the other hand, for all points on the open segment $(Q_2, Q_3)$, the multiplex $\bm{\mathcal{O}}$ defined above is the unique $\bm \theta$-extremal multiplex. Also, at the point $Q_2$ both $\bm{\mathcal{R}}$ and $\bm{\mathcal{O}}$ are $\bm \theta$-extremal. Furthermore, at the point $Q_1$ there are 3 
    $\bm \theta$-extremal multiplexes: the original multiplex $\bm{\mathcal{R}}$, the graph $R^{(1)}$ (which is the blue-colored triangle in layer 1), and the cross layer triangle on the vertices $\{1,2, 3\}$ formed by the blue edges $\{(1, 2), (1, 3)\}$ in layer 1 and the red edge $(2, 3)$ in layer 2.
\end{example}

\begin{example}
    Consider the multiplex shown in Figure \ref{fig:H} (a):
    \begin{align}\label{eq:H}
        \bm{\mathcal{P}} = (\{1, 2, 3, 4, 5\}, \mathcal{P}^{(1)}, \mathcal{P}^{(2)} ) ,
    \end{align}
    where $E(\mathcal{P}^{(1)}) =  \{(1, 2), (2, 3), (2, 5), (3, 4), (4, 5)\}$, which is the graph in layer 1 with edges colored in blue, and $E(\mathcal{P}^{(2)}) =  \{(2, 3), (1, 3) \}$, which is the graph in layer 2 with edges colored in red.  
    \begin{figure}
        \begin{subfigure}[b]{0.18\linewidth}
            \begin{center}
                \begin{tikzpicture}[thick, scale=0.65, every node/.style={draw, circle, fill=gray!30, inner sep=1pt, minimum size=8mm}]




                    \node (1) at (0,1) {1};
                    \node (2) at (-1.5,-1.5) {2};
                    \node (3) at (1.5,-1.5) {3};
                    \node (4) at (1.5, -4) {4};
                    \node (5) at (-1.5, -4) {5};

                    \draw[blue, thick] (1) -- (2);
                    \draw[red, thick] (1) -- (3);
                    \draw[blue, thick] (2) -- (3);
                    \draw[blue, thick] (2) -- (5);
                    \draw[blue, thick] (3) -- (4);
                    \draw[blue, thick] (5) -- (4);
                    \draw[red, thick] (2) to[bend right=40] (3);

                \end{tikzpicture} \\
                { (a) }
            \end{center}
        \end{subfigure}
        \hfill
        \begin{subfigure}[b]{0.45\linewidth}
            \begin{center}
                \includegraphics[width=\linewidth]{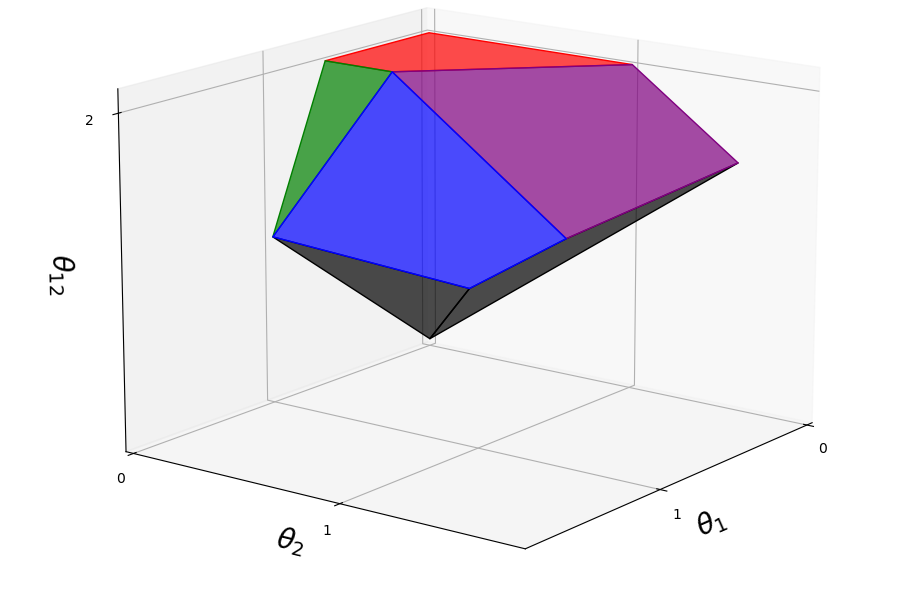} \\
                { (b) }
            \end{center}
        \end{subfigure}
        \hfill
        \begin{subfigure}[b]{0.35\linewidth}
            \begin{center}
                \input{plots/Multiplex_2.tex} \\
                {(c) }
            \end{center}
        \end{subfigure}
        \caption{ \small{(a) The multiplex $\bm{\mathcal{P}}$ defined in \eqref{eq:H}, (b) the full phase diagram in 3D, (c) the 2D phase diagram assuming $\theta_1=\theta_2= \theta$. }}
        \label{fig:H}
    \end{figure}
    The satisfiability region for this multiplex is shown in Figure \ref{fig:H} (b) and the 2-dimensional slice where $\theta_1=\theta_2= \theta$ is shown in Figure \ref{fig:H} (c). In this case, the satisfiability region $\partial_{\bm{\mathcal{P}}}^+$ is the blue polygonal region together with the blue line segment $(O, Q_1)$ (excluding the origin $O$ and the point $Q_1$) as shown in Figure \ref{fig:H} (c). The boundary curve $\partial_{\bm{\mathcal{P}}}$ is the red polyline segment $[Q_1, Q_2, Q_3, Q_4)$ (including the points $Q_1$, $Q_2$, $Q_3$ and excluding the point $Q_4$) and $\partial_{\bm{\mathcal{P}}}^{-}$ is the grey region in Figure \ref{fig:H} (c). In this case, the $\bm \theta$-extremal multiplexes on the boundary are:
    \begin{itemize}
        \item on the open line segment $(Q_1,Q_2)$, $\bm{\mathcal{P}}$
              is uniquely $\bm \theta$-extremal,
        \item on $(Q_2,Q_3)$, the submultiplex induced on the vertices $\{1, 2, 3\}$ is uniquely $\bm \theta$-extremal,
        \item on $(Q_3,Q_4)$, the submultiplex induced on the vertices $\{2, 3\}$ (which is isomorphic to the multiplex $\bm{\mathcal{O}}$ defined in Example \ref{example:edgetriangle}) is uniquely $\bm \theta$-extremal.
    \end{itemize}
    Moreover, at the point $Q_2$ both $\bm{\mathcal{P}}$ and the submultiplex induced on the vertices $\{1, 2, 3\}$ are $\bm \theta$-extremal. Similarly, at $Q_3$ both the submultiplexes induced on the vertices $\{1, 2, 3\}$ and $\{2,3\}$ are $\bm \theta$-extremal. Furthermore, at  $Q_1$ there are three distinct $\bm \theta$-extremal submultiplexes: the original multiplex $\bm{\mathcal{P}}$ and two non-induced submultiplexes $\tilde{\bm{\mathcal{P}}}$ and $\bm{\mathcal{P}}'$ defined as follows:
    \begin{align*}
        \tilde{\bm{\mathcal{P}}} = \left(\{1,2,3,4,5\}, E(\mathcal{P}^{(1)}), \{(1,3)\}\right) \quad \text{and} \quad \bm{\mathcal{P}}' = \left(\{1,2,3,4,5\}, E(\mathcal{P}^{(1)})\setminus\{(2,3)\}, E(\mathcal{P}^{(2)})\right).
    \end{align*}
\end{example}

\begin{remark}[Boundedness of $\partial_{\bm H}^+$] Note that the threshold for the existence of a single edge in the Erd\H{o}s-R\'enyi model $G(n, p_1)$ is $n^2 p_1 \gtrsim 1$. Hence, whenever the multiplex $\bm H$ has an edge that belongs to both the layers, one must have $n^2p_{12} \gtrsim 1$ for $\bm H$ to appear in $G(n, p_1, p_2, p_{12})$. In the parametrization \eqref{eq:exponent}, this means $\theta_{12} \leq 2$ and, by the feasibility constraint, $\max\{\theta_1, \theta_2\} \leq 2$. Hence, the satisfiability region $\Delta_{\bm{H}}$ will be a bounded polyhedral region whenever $\bm{H}$ has an edge that belongs to both layers.  
\end{remark}

\subsection{Asymptotic Normality in the Satisfiability Region }
\label{sec:distributionZHGn}

Once the threshold for the emergence of a multiplex $\bm{H}$ is obtained, the next natural question is to determine the limiting distribution of $X(\bm{H}, \bm G_n)$ in the satisfiability region (where $\bm{G}_n$ contains a diverging number of copies of $\bm {H}$). For the Erd\H{o}s-R\'enyi random graph $G(n, p_1)$ it is well known that the number of copies of a fixed subgraph $H$ is asymptotically normal, whenever $G_n$ contains a diverging number of copies of $H$ (which corresponds to the condition $n p_1^{m_H} \gg 1$) \cite{rucinski1988small,barbour1989central,nowicki1989asymptotic,janson1990functional} (see \cite[Chapter 6]{janson2011random} for a comprehensive discussion).  In this section we will derive the analogous result for submultiplexes, under the parametrization \eqref{eq:exponent}.

To state the result, recall the definition of the Wasserstein distance between two random variables $X$ and $Z$, to be denoted by $d_{\mathrm{Wass}}(X, Z)$, from \cite[Chapter 4]{ChenGoldsteinShao2011}.
%
%
Define
\begin{align}\label{eq:ZHGn}
    Z(\bm H, \bm G_n) := \frac{ X(\bm H, \bm G_n) - \E [X(\bm H, \bm G_n) ]}{ \sqrt{ \mathrm{Var}[X(\bm H, \bm G_n)] } } .
\end{align}
The following theorem establishes the asymptotic normality of $Z(\bm H, \bm G_n)$ with a rate of convergence in terms of the Wasserstein distance.

\begin{theorem}
    \label{thm:ZHGn}
    Suppose $\bm H = (V(\bm{H}), H^{(1)}, H^{(2)})$ is a fixed multiplex with at least one edge. Then, with the parametrization \eqref{eq:exponent}, the following holds: 
    \begin{align}\label{eq:ZHGn}
    d_{\mathrm{Wass}} (Z(\bm H, \bm G_n), N(0, 1) ) \lesssim_{\bm H}  \frac{1}{\sqrt{n^{\Delta_{\bm H}(\bm \theta)}}}  ,  
    \end{align}
    where $\Delta_{\bm H}(\bm \theta)$ is defined in \eqref{eq:parameterH}. This implies that $Z(\bm H, \bm G_n) \stackrel{D} \rightarrow N(0, 1)$, whenever $\Delta_{\bm H}(\bm \theta) > 0$.
\end{theorem}

The proof of  Theorem \ref{thm:ZHGn}  is given in Section \ref{sec:ZHGnpf}. The proof uses Stein's method based on dependency graphs, which bounds the distance to normality in the Wasserstein distance. This result shows that the number of copies of any submultiplex $\bm H$ in the correlated Erd\H{o}s-R\'{e}nyi random multiplex $\bm{G}_n$ is asymptotically normal, for all parameter values for which the count of $\bm H$ in $\bm G_n$ diverges in probability, in the parametrization \eqref{eq:exponent}. This is depicted by the colored polyhedral regions in Figures \ref{fig:edgetriangle} and \ref{fig:H} for their respective examples.

\subsection{ Distribution on the Boundary Surface }
\label{sec:distributionH}

In this section we consider the regime where $\Phi_{\bm H} \asymp 1$ (or in the parametrization \eqref{eq:exponent}, $\Delta_{\bm H} = 0$). This is the boundary between the regimes in which the probability of containing $\bm H$ in $\bm{G}_n$ tends to $0$ and $1$, respectively. Deriving the distribution of $X (\bm H, \bm G_n)$ in this regime is an intriguing problem. To motivate our results, let us begin by recalling the analogous results for the Erd\H{o}s-R\'enyi model $G(n, p_1)$ when $p_1$ is at the threshold (for a particular subgraph $H$).
More precisely, consider a fixed subgraph $H$ and $G_n \sim G(n, p_1)$, such that $n p_1^{m(H)} \rightarrow c \in (0, \infty)$ satisfies the threshold condition (recall the definition of $m(H)$ from \eqref{eq:mH}). Then $X(H, G_n)$ has the following properties:

\begin{enumerate}[label=(\textsf{P\arabic*})]

    \item
          \label{result:balancedH}
          If $H$ is strictly balanced (recall Remark \ref{remark:graphdistribution}), then $X(H, G_n)$ converges to a Poisson distribution (see \cite[Theorem 3.19]{janson2011random}).

    \item
          \label{result:unbalancedH}
          If $H$ is unbalanced, then $\E[X(H, G_n)] \gg 1$ at the threshold. In this case, there is a deterministic sequence $a_n(H) \gg 1$ such that the asymptotic distribution of $\frac{X(H, G_n)}{a_n(H) }$ coincides with that of $X(H', G_n)$, for a certain balanced subgraph $H'$ of $H$. Hence, the problem reduces to the balanced case (see \citep[Section 4]{rucinski1990survey} for details).

    \item If $H$ is balanced, but not strictly balanced, then the limiting distribution of $X(H, G_n)$ can, in principle, be described (using the notion
          of grading \cite{bollobas1989subgraph}), but it has no universal compact form. In particular, it  depends on the structure and the number of subgraphs that maximize \eqref{eq:mH} (see \cite[Section 3.3]{janson2011random} for examples).

\end{enumerate}

Our goal in this section is to derive analogous results for submultiplexes in the correlated Erd\H{o}s-R\'{e}nyi model. The first step towards this is to define the notion of balance in the multiplex setting. Note that since we have a threshold surface, the definition of balance will depend on which parameter value on the surface is chosen. 
Fro this, recall the definition of $\ell_{\bm{F}}(\cdot)$ from \eqref{eq:VF}.

\begin{definition}
    Fix a multiplex $\bm H =(V(\bm H), H^{(1)}, H^{(2)})$ and a point $\bm \theta= (\theta_1, \theta_2, \theta_{12}) \in \partial_{\bm H}$.  Then $\bm H$ is said to be $\bm \theta$-{\it balanced} if $\ell_{\bm H}(\bm \theta) = 0$ and
    \begin{align}\label{eq:balancedH}
        \ell_{\bm F}(\bm \theta) \geq 0 \text{ for every proper submultiplex } \bm F \subset \bm H, \text{ with } |E(\bm F)| \geq 1  .   
    \end{align}
    If the inequality in \eqref{eq:balancedH} is strict, then $\bm H$ is said to be {\it strictly} $\bm \theta$-{\it balanced}. Furthermore, if $\ell_{\bm H}(\bm \theta) > 0$
    and $\ell_{\bm F}(\bm \theta) = 0$, for some submultiplex $\bm F \subset \bm H$, with $|E(\bm F)| \geq 1$, 
    then $\bm H$ is said to be $\bm \theta$-{\it unbalanced}.
\end{definition}

\begin{remark}Note that $\bm H$ is $\bm \theta$-balanced if and only if $\bm H$ is $\bm \theta$-extremal (recall Definition \ref{defn:boundaryH}).
\end{remark}

We begin with the case when $\bm{H}$ is strictly $\bm \theta$-balanced. Let $\mathrm{Aut}(\bm H)$ denote the collection of all permutations of $V(\bm H)$ that
preserve the edge sets in both layers. 
Then we have the following result, which is proved in Section \ref{sec:balanceHpf}.

\begin{theorem}
    \label{thm:balancedH}
Fix a multiplex $\bm H =(V(\bm H), H^{(1)}, H^{(2)})$ and a point $\bm \theta= (\theta_1, \theta_2, \theta_{12}) \in \partial_{\bm H}$ such that $\bm{H}$ is strictly $\bm \theta$-balanced and $\theta_{12} > \max\{\theta_1, \theta_2\}$. 
    Then
    \begin{align*}
        X(\bm{H}, \bm{G}_n) \xrightarrow{D} \mathrm{Pois}\left(\frac{1}{|\mathrm{Aut}(\bm H)|}\right).
    \end{align*} 
\end{theorem}

\begin{remark}  
The condition $\theta_{12}>\max\{\theta_1,\theta_2\}$ cannot, in general, be omitted if one wishes to obtain a Poisson limit. To see this, consider the multiplex
$\bm H = (\{1,2,3\}, H^{(1)}, H^{(2)} )$,  where
$E(H^{(1)})=\left\{ (1,2) \right\}$ and $E(H^{(2)})=\left\{ (1,3) \right\}$.
Thus, $\bm H$ consists of a two-edge path whose two edges belong
to different layers. Suppose that $\theta_1=\theta_2=\theta_{12}=\frac{3}{2}.$ In this case, $\bm H$ is strictly $\bm\theta$-balanced and the two layers of $\bm G_n$ coincide almost surely. Let $Y_n$ denote the number of copies of the 
two-edge path $P_3$ in the common graph $G_n^{(1)}=G_n^{(2)} \sim G(n,n^{-\frac{3}{2}})$.  Each copy of $P_3$ gives rise to exactly two distinct copies of
$\bm H$ in $\bm G_n$, corresponding to the two possible
assignments of its two edges to the two layers. Therefore,
$X(\bm H,\bm G_n)=2Y_n$.  Since $P_3$ is strictly balanced and $\mathbb E[Y_n] \rightarrow\frac{1}{2}$,  it follows from \cite[Theorem~3.19]{janson2011random} that $Y_n\xrightarrow{D}\operatorname{Pois}(\frac{1}{2})$. Hence, 
$X(\bm H,\bm G_n) \stackrel{D} \rightarrow 2 \mathrm{Pois}(\frac{1}{2})$, which is a scaled Poisson distribution.  
\end{remark}

Next, we consider the case where $\bm H$ is $\bm{\theta}$-unbalanced. In this case, we can rescale $X(\bm{H}, \bm G_n)$ such that its asymptotic distribution coincides with that of a $\bm{\theta}$-balanced submultiplex of $\bm{H}$. Towards this we need a few definitions: To begin with, fix a point $\bm \theta= (\theta_1, \theta_2, \theta_{12}) \in \partial_{\bm H}$ and recall the notion of $\bm \theta$-extremal from Definition \ref{defn:boundaryH}. Denote the collection of  $\bm \theta$-extremal submultiplexes of $\bm H$ as $\mathcal{A}_{\bm H}(\bm \theta)$. The {\it core} of $\bm H$ is defined as the largest $\bm \theta$-extremal submultiplex of $\bm H$, that is, the unique maximum element in $\mathcal{A}_{\bm H}(\bm \theta)$ in terms of multiplex containment (which is well defined by Lemma \ref{lem:extreme}).\footnote{For two multiplexes  $\bm{G} = (V(\bm{G}), G^{(1)}, G^{(2)})$ and
    $\bm{H} = (V(\bm{H}), H^{(1)}, H^{(2)})$, their union
    $\bm{G} \cup \bm{H}$ is the multiplex with vertex set $V(\bm{G}) \cup V(\bm{H})$ and two layers with edge sets  $E(G^{(1)}) \cup E(H^{(1)}) $ and $E(G^{(2)}) \cup E(H^{(2)})$, respectively. The intersection $\bm{G} \cap \bm{H}$ is defined similarly (with the unions replaced by intersections).}  We also require the notion of an extension of a multiplex. Let
$\bm F\subseteq\bm H$ be a fixed multiplex, and let
$\bm F'$ be a copy of $\bm F$. An \emph{extension} of
$\bm F'$ to $\bm H$ is a copy $\bm H'$ of
$\bm H$ such that $\bm F'\subseteq\bm H'$.  Given a fixed copy $\overline{\bm H}'$ of the core
$\overline{\bm H}$, let $\ext_n(\overline{\bm H}';\overline{\bm H},\bm H)$
denote the number of distinct copies $\bm H'$ of
$\bm H$ satisfying $\overline{\bm H}'\subseteq\bm H'
\subseteq
\bm G_n\cup\overline{\bm H}'$. By exchangeability, the
distribution of this random variable does not depend on the
particular copy $\overline{\bm H}'$. We therefore write $\ext_n(\overline{\bm H},\bm H)$ for a random variable having this common distribution.

\begin{theorem}
Fix a multiplex $\bm H =(V(\bm H), H^{(1)}, H^{(2)})$ and a point $\bm \theta= (\theta_1, \theta_2, \theta_{12}) \in \partial_{\bm H}$. Then, with the parametrization \eqref{eq:exponent}, 
    \begin{align}\label{eq:unbalancedH}
        \left|\frac{X(\bm H, \bm G_n)}{\E[\ext_n(\overline{\bm{H}}, \bm{H})]} - X(\overline{\bm H}, \bm G_n) \right| \xrightarrow{P} 0 . 
    \end{align}  
    \label{thm:unbalancedH} 
\end{theorem}

The proof of Theorem \ref{thm:unbalancedH} is given in Section \ref{sec:unbalancedHpf}.
A main technical step in the proof is to show that the core $\overline{\bm H}$ is well-defined. The challenge here is that, unlike in the case of graphs (single-layer networks), the notions of balanced/extremal depend on the parameter value $\bm \theta \in \partial_{\bm H}$. The key observation here is that $\ell_{\bm F}(\bm \theta)$ is modular
(see \eqref{eq:modular} for the definition), for $\bm \theta \in \partial_{\bm H}$, as $\bm F$ ranges over the set of $\bm \theta$-extremal multiplexes. Using this, we can establish that the core is well-defined. Then by estimating the number of extensions of a planted core together with a second-moment argument, the result in \eqref{eq:unbalancedH} follows.

\begin{remark}
    Note that Theorem \ref{thm:unbalancedH} holds irrespective of whether $\bm{H}$ is $\bm{\theta}$-balanced or $\bm{\theta}$-unbalanced.
    \begin{itemize}
        \item When $\bm{H}$ is $\bm{\theta}$-balanced, then $\overline{\bm{H}} = \bm{H}$ and $\ext_n(\overline{\bm{H}}', \bm{H}) = \ext_n(\bm{H}', \bm{H}) = 1$, for any copy $\bm{H}'$ of $\bm{H}$. Hence, in this case, \eqref{eq:unbalancedH} holds trivially.
        \item When $\bm{H}$ is $\bm{\theta}$-unbalanced, then $\overline{\bm{H}} \subset \bm{H}$ and $\E[X(\bm H, \bm G_n)] \gg 1$. In this case,  \eqref{eq:unbalancedH}  shows that the asymptotic distribution of $\frac{X(\bm H, \bm G_n)}{\E[\ext_n(\overline{\bm{H}}, \bm{H})]}$ is the same as that of $X(\overline{\bm H}, \bm G_n)$. Since $\bm \theta \in \partial_{\bm H}$ and $\overline{\bm{H}}$ is $\bm{\theta}$-extremal, we also have $\bm \theta \in \partial_{\overline{\bm H}}$ and $\overline{\bm{H}}$ is $\bm{\theta}$-balanced. Hence, \eqref{eq:unbalancedH} reduces the problem to the balanced case.
    \end{itemize}
\end{remark}

A useful special case is when $\bm H =(V(\bm H), H^{(1)}, H^{(2)})$ is $\bm{\theta}$-unbalanced such that $|\mathcal A_{\bm H}(\bm \theta)| = 1$, that is, there is a unique $\bm \theta$-extremal submultiplex (which is not $\bm H$ itself). Then the core $\overline{\bm{H}}$ is strictly $\bm{\theta}$-balanced and invoking Theorems \ref{thm:balancedH} and \ref{thm:unbalancedH} we have the following result:

\begin{cor}\label{cor:uniqueH}
    Fix a multiplex $\bm H =(V(\bm H), H^{(1)}, H^{(2)})$ and a point $\bm \theta= (\theta_1, \theta_2, \theta_{12}) \in \partial_{\bm H}$ such that $\theta_{12} > \max\{\theta_1, \theta_2\}$ and $|\mathcal A_{\bm H}(\bm \theta)| = 1$. Then, 
 \begin{align*}
        \frac{X(\bm H, \bm G_n)}{\E[\ext_n(\overline{\bm{H}}, \bm{H})]} \xrightarrow{D} \mathrm{Pois}\left(\frac{1}{ |\mathrm{Aut}(\overline{\bm H})| }\right) .
    \end{align*}
\end{cor}


To illustrate the above results, let us return to the edge-triangle multiplex from Example \ref{example:edgetriangle}.

\begin{example}[Example \ref{example:edgetriangle} continued]
    Recall the edge-triangle multiplex $\bm{\mathcal{R}}$ from \eqref{eq:H123}. Assume $\theta_1=\theta_2=\theta$. Note that $|\mathrm{Aut}(\bm{\mathcal{R}})|=2$. In this case, the threshold curve $\partial_{\bm{\mathcal{R}}}$ is the red polyline segment $[Q_1, Q_2, Q_3)$ (including the point $Q_1$ and excluding the point $Q_3$) as shown in Figure \ref{fig:edgetriangle} (c). Depending on the location of the point $\bm \theta = (\theta, \theta_{12})$ on this polyline segment we have the following cases:

    \begin{itemize}

        \item $\bm \theta$ lies on the open segment $(Q_1, Q_2)$. Here, $\bm{\mathcal{R}}$ is the unique $\bm \theta$-extremal multiplex. Hence, on this line segment $\bm{\mathcal{R}}$ is strictly $\bm \theta$-balanced and Theorem \ref{thm:balancedH} implies:
              $$X(\bm{\mathcal{R}}, \bm G_n) \stackrel{D} \rightarrow \mathrm{Pois}(\tfrac{1}{2}).$$

        \item $\bm \theta$ lies on the open segment $(Q_2, Q_3) = \{(\theta, 2) : 0 < \theta < \tfrac{1}{2} \}$. Here, the multiplex $\bm{\mathcal{O}}$ defined in Example \ref{example:edgetriangle} is the unique $\bm \theta$-extremal submultiplex. Hence, on this line segment $\bm{\mathcal{R}}$ is $\bm \theta$-unbalanced with $|\mathcal A_{\bm{\mathcal{R}}}(\bm \theta)| = 1$, which means Corollary \ref{cor:uniqueH} can be applied. For this, we need to compute $\E[\ext_n(\overline{\bm{\mathcal{R}}}, \bm{\mathcal{R}})]$. To this end, without loss of generality assume that a copy of the multiplex $\bm{\mathcal{O}}$ is planted on the vertices $\{1, 2\}$ of $\bm G_n$. To extend this copy of $\bm{\mathcal{O}}$ to a copy of $\bm{\mathcal{R}}$ we need to choose another vertex $v \in [n]\backslash\{1, 2\}$ (for which we have $(n-2)$ choices) and ensure that the edges $(1, v)$ and $(2, v)$ are present in the first layer of $\bm{G}_n$ (which has probability $p_{1}^2 = n^{-2\theta}$). Hence,   $\E[\ext_n(\overline{\bm{\mathcal{R}}}, \bm{\mathcal{R}})]  \sim n^{1-2 \theta}$. Hence, applying Corollary \ref{cor:uniqueH} gives, 
              \begin{align*}
                  \frac{X(\bm{\mathcal{R}}, \bm G_n)}{n^{1-2 \theta}} \xrightarrow{D} \mathrm{Pois}\left(\tfrac{1}{ 2 }\right) .
              \end{align*}

        \item Next, suppose $\bm \theta = Q_2 = (\frac{1}{2}, 2)$. In this case, both $\bm{\mathcal{R}}$ and $\bm{\mathcal{O}}$ are $\bm \theta$-extremal. This means at this point $\bm{\mathcal{R}}$ is $\bm \theta$-balanced, but not strictly $\bm \theta$-balanced. This setting is not covered by the general results discussed above. In fact, in this case the asymptotic distribution of $X(\bm{\mathcal{R}}, \bm G_n)$ is not Poisson. By a direct calculation it can be shown that
              \begin{align*}
                  X(\bm{\mathcal{R}}, \bm G_n)  \xrightarrow{D}  \sum_{i=1}^Z Z_i ,   \end{align*}
              where $Z \sim \mathrm{Pois}(\frac{1}{2})$ and independently $Z_1, Z_2, \ldots $ are i.i.d. $\mathrm{Pois}(1)$ random variables. The idea behind this is as follows: Note that $\theta_{12} =  2$ implies $p_{12} = \frac{1}{n^2}$. This means the asymptotic distribution of $X(\bm{\mathcal{O}}, \bm{G}_n)$ is $\mathrm{Pois}(\frac{1}{2})$. Moreover, each copy of the multiplex $\bm{\mathcal{O}}$ in $\bm{G}_n$ has asymptotically $\mathrm{Pois}(1)$ number of common neighbors in layer 1 of $\bm G_n$ (since $\theta = \frac{1}{2}$ means $p_1 = \frac{1}{\sqrt n}$), each of which forms a copy of $\bm{\mathcal{R}}$ in $\bm{G}_n$.

        \item The remaining case is $\bm \theta = Q_1 = (1, 1)$. In this case, recall from Example \ref{example:edgetriangle} that there are three $\bm \theta$-extremal submultiplexes (one of which is $\bm{\mathcal{R}}$ itself). This means, at this point $\bm{\mathcal{R}}$ is $\bm \theta$-balanced, but not strictly $\bm \theta$-balanced, as in the previous case. Note that this corresponds to the parameter values $p_1=p_2= p_{12} = \frac{1}{n}$, that is, the layers of $\bm G_n$ are perfectly correlated. Hence, given a triangle in layer 1 of $\bm G_n$, there are 3 ways to choose the red edge in layer 2 of $\bm G_n$ (refer to Figure \ref{fig:edgetriangle} (a)) to obtain a copy of
              $\bm{\mathcal{R}}$ in $\bm G_n$. This implies,
              $$X(\bm{\mathcal{R}}, G_n) = 3 X(R^{(1)}, G_n^{(1)}) \stackrel{D} \rightarrow 3 \mathrm{Pois}(\tfrac{1}{6}) , $$
              since the number of triangles in $G(n, p_1)$ when $np_1 \rightarrow 1$, converges to $\mathrm{Pois}(\frac{1}{6})$ (see \cite[Theorem 3.19]{janson1990functional}).

    \end{itemize}  
\end{example}

\section{Proof of Theorem \ref{thm:threshold}}
\label{sec:thresholdpf}

Let $\bm G_n = ([n], G_n^{(1)}, G_n^{(2)}) \sim G(n, p_1, p_2, p_{12})$ be a sample from the correlated Erd\H{o}s-R\'enyi model as in Definition \ref{definition:correlated}.   Fix a multiplex $\bm H = (V(\bm{H}), H^{(1)}, H^{(2)})$ with at least one edge. Recalling \eqref{eq:NHGnM}, note that if $\mathcal{H}_n=\{\bm H_1,\,\bm H_2,\,\ldots,\, \bm H_M\}$ is the collection of all copies of $\bm H$ in the complete multiplex $\bm{K}_n$ with vertex set $[n]$, where $M = X(\bm{H}, \bm{K}_n)$, then     
\begin{align}\label{eq:XHGnM}
    X(\bm{H}, \bm{G}_n) = \sum_{s=1}^M \bm{1}\{ \bm H_s \subseteq \bm{G}_n \} .
\end{align} 
With this representation, in the following lemma we compute the asymptotic orders of the mean and the variance of $X(\bm{H}, \bm{G}_n)$. For this, it is convenient to define, for a 
multiplex $\bm H = (V(\bm H), H^{(1)}, H^{(2)})$, 
$$\lambda_{p_1, p_2, p_{12}}(\bm H) := p_1^{|E(H^{(1)})\setminus E(H^{(2)})|} p_2^{|E(H^{(2)})\setminus E(H^{(1)})|} p_{12}^{|E(H^{(1)})\cap E(H^{(2)})|}. $$


\begin{lemma}
    \label{lm:meanvariance}
    For any multiplex $\bm H = (V(\bm{H}), H^{(1)}, H^{(2)})$ with at least one edge,
    \begin{align}\label{eq:XHGnexpectation}
        \E \left[ X(\bm{H}, \bm{G}_n) \right]  \asymp_{\bm H} n^{|V(\bm{H})|} \lambda_{p_1, p_2, p_{12}}(\bm H) .
    \end{align}
    Furthermore,
    \begin{align}\label{eq:XHGnvariance}
        \Var \left[ X(\bm{H}, \bm{G}_n) \right] \lesssim_{\bm H} n^{2|V(\bm{H})|} \lambda_{p_1, p_2, p_{12}}(\bm H)^2 \sum_{ \substack{\bm F \subseteq \bm H \\ |E(\bm{F})| \geq 1 }  } \frac{1}{n^{|V(\bm F)| } \lambda_{p_1, p_2, p_{12}}(\bm F)}  .
    \end{align} 
    Further, under the parametrization \eqref{eq:exponent} one has the following: 
        \begin{align}\label{eq:Hvariance}
        \Var \left[ X(\bm{H}, \bm{G}_n) \right] \asymp_{\bm H} n^{2|V(\bm{H})|} \lambda_{p_1, p_2, p_{12}}(\bm H)^2 \sum_{ \substack{\bm F \subseteq \bm H \\ |E(\bm{F})| \geq 1 }  } \frac{1}{n^{|V(\bm F)| } \lambda_{p_1, p_2, p_{12}}(\bm F)}  .
    \end{align} 
\end{lemma}

\begin{proof}
    From \eqref{eq:XHGnM} we have,
    \begin{align}
        \E\left[ X(\bm{H}, \bm{G}_n) \right]  =\sum_{ s= 1 }^M \P( \bm H_s \subseteq \bm{G}_n )       & = X(\bm{H}, \bm{K}_n) p_1^{|E(H^{(1)})\setminus E(H^{(2)}) |}p_2^{|E(H^{(2)})\setminus E(H^{(1)}) |}p_{12}^{|E(H^{(1)})\cap E(H^{(2)}) |}          \nonumber           \\
                                             & \asymp_{\bm H} n^{|V(\bm H)|} \lambda_{p_1, p_2, p_{12}}(\bm H)  ,  
  \label{eq:XHKnp}                  
    \end{align}
    since $X(\bm{H}, \bm{K}_n) \asymp_{\bm H} n^{|V(\bm H)|}$.

    Next, we compute $\Var \left[X(\bm{H}, \bm{G}_n) \right]$. From \eqref{eq:XHGnM},
    \begin{align}
        \Var\left[X(\bm{H}, \bm{G}_n)\right]  & = \sum_{\substack{\bm{H}_s, \bm{H}_t \in \mathcal{H}_n  \\ E(\bm{H}_s) \cap E(\bm{H}_t) \ne \emptyset }}\Cov[ \bm{1}\{\bm{H}_s \subseteq \bm{G}_n  \} , \bm{1}\{\bm{H}_t \subseteq \bm{G}_n  \}  ] \nonumber \\ 
         & =  \sum_{\substack{\bm{H}_s, \bm{H}_t \in \mathcal{H}_n \\ E(\bm{H}_s) \cap E(\bm{H}_t) \ne \emptyset }} \left\{ \P( \bm{H}_s \subseteq \bm{G}_n \text{  and }   \bm{H}_t \subseteq \bm{G}_n ) -  \P( \bm{H}_s \subseteq \bm{G}_n )^2 \right\}.
        \label{var_upbd}
    \end{align}
For $\bm H_s, \bm H_t \in \mathcal H_n$, assume $\bm{H}_s=(V(\bm H_s), H_s^{(1)},H_s^{(2)})$ and $\bm{H}_t=(V(\bm H_t), H_t^{(1)},H_t^{(2)})$. Let $\bm{J}=\bm{H}_s\cup \bm{H}_t$. Then,
     \begin{align}\label{eq:HstGn}
        \P( \bm{H}_s \subseteq \bm{G}_n \text{  and }   \bm{H}_t \subseteq \bm{G}_n ) 
        = \lambda_{p_1, p_2, p_{12}}(\bm J)  .  
    \end{align}
 Now, for notational convenience let $A_s:=E(H_s^{(1)})$, $B_s:=E(H_s^{(2)})$, $A_t:=E(H_t^{(1)})$, $B_t:=E(H_t^{(2)})$.
Then the (one-sided) overlap multiplex $\bm F_{s, t}= (V(\bm F_{s, t}), F_{s, t}^{(1)},F_{s, t}^{(2)})$ is defined as follows:  
$V(\bm F_{s, t}):=V(\bm H_s)\cap V(\bm H_t)$, 
and $ E(F_{s, t}^{(1)}) :=  A_s\cap(A_t\cup B_t) \text{ and }  E(F_{s, t}^{(2)}) :=
B_s\cap(A_t\cup B_t)$.  Thus, $\bm F_{s, t}$ retains every underlying edge of
$\bm H_s$ that is also used by $\bm H_t$, while preserving
its original layer type in $\bm H_s$. In particular, $\bm F_{s, t} \subseteq \bm H_s$,
and, hence it is isomorphic to a submultiplex of $\bm H$. (Note that, since $\bm F_{s, t}$ is one-sided, $\bm F_{s, t} \ne \bm F_{t, s} $, in general, although the following argument would work with $F_{t, s} $ as well.) We now claim that 
    \begin{align}\label{eq:GnHst}
        \P( \bm{H}_s \subseteq \bm{G}_n \text{  and }   \bm{H}_t \subseteq \bm{G}_n ) 
        \leq & \frac{\lambda_{p_1, p_2, p_{12}}(\bm H)^2}{\lambda_{p_1, p_2, p_{12}}(\bm F_{s, t})}  .  
    \end{align} 
    To prove this, define $D_{s,t}^{(1)}:=B_s\cap(A_t\setminus B_t)$, 
and $D_{s,t}^{(2)}:=A_s\cap(B_t\setminus A_t)$. Then, a set-theoretic calculation gives,  
\begin{align*}
\frac{\P(\bm H_s\subseteq\bm G_n \text{ and } \bm H_t\subseteq\bm G_n)}{  \frac{\lambda_{p_1, p_2, p_{12}}(\bm H)^2}{\lambda_{p_1, p_2, p_{12}}(\bm F_{s, t})}  }   & = \left(\frac{p_{12}}{p_1}\right)^{|D_{s,t}^{(1)}|}
\left(\frac{p_{12}}{p_2}\right)^{|D_{s,t}^{(2)}|}  
\le 1  ,    
\end{align*}
where the last step uses $p_{12} \leq \min \{p_1, p_2\}$. Combining the above proves \eqref{eq:GnHst}.

Furthermore, $\mathbb P\left(\bm H_s\subseteq\bm G_n\right)^2=
\lambda_{p_1,p_2,p_{12}}(\bm H)^2$. Therefore,
\begin{align}
\mathrm{Cov}\left[ \bm 1\{\bm H_s\subseteq\bm G_n\}, \bm 1\{\bm H_t\subseteq\bm G_n\} \right] \leq \mathbb P\left( \bm H_s\subseteq\bm G_n \text{ and } \bm H_t\subseteq\bm G_n \right) &\leq \frac{\lambda_{p_1,p_2,p_{12}}(\bm H)^2}{\lambda_{p_1,p_2,p_{12}}(\bm F_{s,t})}.  
\label{eq:covarianceXHGn}
\end{align}
For every fixed submultiplex $\bm F\subseteq\bm H$, the number
of ordered pairs $(\bm H_s,\bm H_t)$ such that $\bm F_{s, t}\cong\bm F$ is at most $\lesssim_{\bm H} n^{2|V(\bm H)|-|V(\bm F)|}$. Combining this with \eqref{var_upbd} and \eqref{eq:covarianceXHGn}, the result in \eqref{eq:XHGnvariance} follows.

To prove \eqref{eq:Hvariance} we have to show the corresponding lower bound also holds, under the parametrization \eqref{eq:exponent}. To begin with, for $U \subset V(\bm H)$, denote by $\bm H[U]$ the submultiplex of $\bm H$ induced on the set of vertices $U$. Then 
$$\sum_{\substack{U\subseteq V(\bm H)\\ |E(\bm H[U])|\geq1}} \frac{1}{n^{|U|}\lambda_{p_1, p_2, p_{12}}(\bm H[U])} \leq \sum_{\substack{\bm F\subseteq\bm H\\|E(\bm F)|\geq1}} \frac{1}{n^{|V(\bm F)|}\lambda_{p_1, p_2, p_{12}}(\bm F)}.$$
Further, using $p_{12}\leq\min\{p_1,p_2\}\leq1$ gives, $ n^{|V(\bm F)|}\lambda_{p_1, p_2, p_{12}}(\bm F) \geq n^{|V(\bm F)|} \lambda_{p_1, p_2, p_{12}}(\bm H[V(\bm F)])$.  
Since $\bm H$ is fixed, there are only $\lesssim_{\bm H}  1$ submultiplexes with a prescribed vertex set. Hence, 
\begin{equation}\label{eq:induced-sum-equivalence} 
\sum_{\substack{\bm F\subseteq\bm H\\|E(\bm F)|\geq1}} \frac{1}{n^{|V(\bm F)|}\lambda_{p_1, p_2, p_{12}}(\bm F)}  \asymp_{\bm H} \sum_{\substack{U\subseteq V(\bm H)\\ |E(\bm H[U])|\geq1}} \frac{1}{n^{|U|}\lambda_{p_1, p_2, p_{12}}(\bm H[U])}  .  
\end{equation}
Now, fix $U\subseteq V(\bm H)$ and denote by $\tilde{\mathcal P}_U$, the collection of ordered pairs of injective maps $\varphi_s,\varphi_t:V(\bm H)\longrightarrow [n]$, 
satisfying $\varphi_s|_U=\varphi_t|_U$ and $\varphi_s(V(\bm H))\cap\varphi_t(V(\bm H)) =\varphi_s(U)=\varphi_t(U)$.
Define $\mathcal P_U$ to be the collection of ordered pairs of copies $\left(\bm H_{\varphi_s},\bm H_{\varphi_t}\right)$
arising from pairs $(\varphi_s, \varphi_t)\in \tilde{\mathcal P}_U$.
Note that for every $(\bm H_s, \bm H_t) \in\mathcal P_U$, the common edges of $\bm H_s$ and $\bm H_t$ are precisely the edges of $\bm H[U]$, with the same layer types in the two copies. Consequently,
\begin{align*}
\Cov[\bm 1\{\bm H_s\subseteq\bm G_n\}, \bm 1\{\bm H_t\subseteq\bm G_n\}] &=
\frac{\lambda_{p_1, p_2, p_{12}}(\bm H)^2}{\lambda_{p_1, p_2, p_{12}}(\bm H[U])}
(1-\lambda_{p_1, p_2, p_{12}}(\bm H[U])).
\end{align*}
Since $\bm H[U]$ contains an edge, under the scaling \eqref{eq:exponent} it follows that $\lambda_{p_1, p_2, p_{12}}(\bm H[U]) \ll 1$. Combining this with the bound 
$|\mathcal P_U| \asymp_{\bm H} n^{2 |V(\bm H)|-|U|}$ gives,  
\begin{align}
\sum_{\substack{U\subseteq V(\bm H)\\|E(\bm H[U])|\geq1}}
\ \sum_{ (\bm H_s, \bm H_t) \in\mathcal P_U}
& \Cov[\bm 1\{\bm H_s\subseteq\bm G_n\}, \bm 1\{\bm H_t\subseteq\bm G_n\}]  \nonumber  \\ 
& \gtrsim_{\bm H}
n^{2 |V(\bm H)|} \lambda_{p_1, p_2, p_{12}}(\bm H)^2   \sum_{\substack{U\subseteq V(\bm H)\\ |E(\bm H[U])|\geq1}} \frac{1}{n^{|U|}\lambda_{p_1, p_2, p_{12}}(\bm H[U])} .
\label{eq:positive-covariance-contribution}
\end{align} 
Note that a fixed ordered pair of copies can belong to $\mathcal P_U$ for at most $\lesssim_{\bm H} 1$ choices of $U$. Consequently, after dividing by this bounded multiplicity, if necessary, the LHS of \eqref{eq:positive-covariance-contribution} is bounded above by a constant multiple of the sum of the positive covariance terms in the variance expansion. Further, the contribution from the negative covariance terms is negligible relative to
\eqref{eq:positive-covariance-contribution} under \eqref{eq:exponent}. Hence, from \eqref{eq:induced-sum-equivalence} and \eqref{eq:positive-covariance-contribution}, the required lower bound in \eqref{eq:Hvariance} follows. 
\end{proof}

\color{black}

With the above lemma, we can complete the proof of Theorem \ref{thm:threshold}. First, suppose $\Phi_{\bm H} \ll 1$. Let $\bm{F}$ be a submultiplex where the minimum for $\Phi_{\bm H}$ is attained (recall \eqref{eq:PhiH}). Then, from \eqref{eq:XHGnexpectation},
\begin{align*}
    \P(X(\bm{H}, \bm{G}_n)>0)  \leq \P(X(\bm{F}, \bm{G}_n)>0) \leq \E[X(\bm{F},\bm{G}_n)] \asymp_{\bm H} \Phi_{\bm{H}} \ll 1 .
\end{align*}
This proves the 0-statement in Theorem \ref{thm:threshold}.

Next, suppose $\Phi_{\bm{H}} \gg 1$.
%
%
Then by Lemma \ref{lm:meanvariance},
\begin{align}\label{eq:XHGn2moment}
    \frac{\Var \left[X(\bm{H}, \bm{G}_n) \right] }{(\E[X(\bm{H}, \bm{G}_n)])^2} &  \sum_{ \substack{\bm F \subseteq \bm H \\ |E(\bm{F})| \geq 1 }  } \frac{1}{n^{|V(\bm F)| } \lambda_{p_1, p_2, p_{12}}(\bm F)}  \nonumber \\  
                                                                                & \lesssim_{\bm H} \frac{1}{ \min_{ \substack{\bm F \subseteq \bm H                       \\ |E(\bm{F})| \geq 1 } } n^{|V(\bm{F})|} \lambda_{p_1, p_2, p_{12}}(\bm F)} = \frac{1}{ \Phi_{\bm{H}} } \ll 1. 
\end{align}
This means, $\frac{X(\bm H, \bm G_n)}{\E[X(\bm H, \bm G_n)]} \stackrel{P} \rightarrow 1$, which implies $X(\bm{H}, \bm{G}_n) \stackrel{P} \rightarrow \infty$,
since $\E[X(\bm H,G_n)] \rightarrow \infty$, whenever $\Phi_{\bm H} \gg 1$. This completes the proof of Theorem \ref{thm:threshold}.  \hfill $\Box$

\section{ Proof of Theorem \ref{thm:ZHGn} }
\label{sec:ZHGnpf}

The proof of Theorem \ref{thm:ZHGn} will use Stein's method based on  dependency graphs. We begin by defining the notion of a dependency graph.

\begin{definition}[Dependency graph]
    Let $\{X_v\}_{v\in V}$ be a family of random variables (defined on some common probability space) indexed by a finite set $V$. Then a graph $\mathcal G$ with vertex set $V$ is said to be a dependency graph for the collection $\{X_v\}$ if the following holds: For any two subsets of vertices $A, B \subseteq V$ such that there is no edge in $\mathcal{G}$ from any vertex in $A$ to any vertex in $B$, the collections of random variables $\{X_v\}_{v\in A}$ and $\{X_v\}_{v\in B}$ are independent.
\end{definition}

For a dependency graph $\mathcal G$ and $u \in V$, denote by $\overline{N}_\mathcal{G}(u)$ the set of neighbors of $u$ in $\mathcal{G}$ and $u$ itself. Also,  denote $\overline{N}_\mathcal{G}(u,\,v):=\overline{N}_\mathcal{G}(u)\cup\overline{N}_\mathcal{G}(v).$ We will use the following version of Stein's method based on dependency graph.

\begin{theorem}[{\cite[Theorem 6.33]{janson2011random}}]
    \label{thm:W}
    Suppose $\{X_v\}_{v\in V}$ be a family of random variables with dependency graph $\mathcal{G}$. Assume $\E[X_v]=0$, for all $v\in V.$ Let $W=\frac{1}{\sigma}\sum_{v\in V}X_v$,
    where $\sigma^2=\Var[\sum_{v\in V}X_v]$.
    Assume, for all $u,\,v\in V,$ the following bounds hold:   
    \begin{equation*}
        \sum_{v\in V}\E[|X_v|] \le R \quad\text{and}\quad\sum_{w\in\overline{N}_\mathcal{G}(u,\,v)}\E[|X_w|\mid X_u,\,X_v] \le Q .
    \end{equation*}
    Then, for $Z \sim N(0, 1)$,  $d_\mathrm{Wass}(W, Z)\lesssim \frac{RQ^2}{\sigma^3}$.
\end{theorem}

With the above preparations we can now proceed with the proof of Theorem \ref{thm:ZHGn}. As before, let $\mathcal{H}_n=\{\bm H_1,\,\bm H_2,\,\ldots,\, \bm H_M\}$ be the collection of all copies of $\bm H$ in $\bm{K}_n$, where $M = X(\bm{H}, \bm{K}_n)$.  
For $1\le s \le M$, define $I_s:=\bm{1}\{ \bm H_s \subseteq \bm G_n\}$ and
\begin{align*}
    X_s =I_s-\mathbb{E}[I_s]=I_s - \mu ,
\end{align*}
where $\mu = \lambda_{p_1, p_2, p_{12}}(\bm H)$ (recall \eqref{eq:XHKnp}). Define $\sigma^2= \Var[X(\bm H,\, \bm G_n)]$. Recalling \eqref{eq:ZHGn}, note that
$Z(\bm H, \bm G_n) = \frac{1}{\sigma} \sum_{s=1}^{M}X_s$. Construct a dependency graph $\mathcal G$ of the collection of random variables $\{X_s : \bm H_s \in\mathcal{H}_n\}$ on the vertex set $\{1,\,2,\,\ldots,\,M\}$ as follows:  Connect the edge $(s,\,t)$, for $1 \leq s < t \leq M$ in $\mathcal{G}$ whenever $|E(\bm H_s)\cap E(\bm H_t)| \geq 1$. (Recall that for any multiplex $\bm H = (V(\bm H), H^{(1)}, H^{(2)})$, $E(\bm H):=E(H^{(1)})\cup E(H^{(2)}$.)

To apply Theorem \ref{thm:W}, first note that $\E|X_s| = 2\mu(1-\mu) \leq2 \mu$. Hence, $R_n := \sum_{s=1}^M\E|X_s| \lesssim_{\bm H} M\mu = \E[X(\bm H,\bm G_n)],$ where the last step follows from \eqref{eq:XHKnp}.  
Next, fix $s, t \in \{1, 2, \ldots, M\}$. For $w\in\overline N_{\mathcal G}(s,t)$, define the (one-sided) overlap multiplex $\bm H_{w, \{s,t\} }= (V(\bm H_{w, \{s,t\} }), H_{w, \{s,t\} }^{(1)},
H_{w, \{s,t\} }^{(2)} )$ as: $V(\bm H_{w, \{s,t\} }):= V(\bm H_w)\cap (V(\bm H_s)\cup V(\bm H_t) )$ and
$$E(H_{w, \{s,t\} }^{(1)}) :=
E(H_w^{(1)})\cap ( E(\bm H_s)\cup E(\bm H_t) ) \text{ and } E(H_{w, \{s,t\} }^{(2)})
:=E(H_w^{(2)})\cap ( E(\bm H_s)\cup E(\bm H_t) ) . $$
Note that $\bm H_{w, \{s,t\} }\subseteq\bm H_w$. Further, since
$w\in\overline N_{\mathcal G}(s,t)$, $|E(\bm H_{w, \{s,t\} })|\geq1$. Let $\mathcal A_{w, \{s,t\}}$ denote the event that all edge
requirements of $\bm H_w$ whose underlying vertex pairs lie
outside $E(\bm H_s)\cup E(\bm H_t)$ are satisfied. Then $\{\bm H_w\subseteq\bm G_n\}
\subseteq\mathcal A_{w, \{s,t\}}$. Note that the event $\mathcal A_{w, \{s,t\}}$ depends only on the edge vectors outside $E(\bm H_s)\cup E(\bm H_t)$. Hence, 
\begin{align*}
\mathbb P(\bm H_w\subseteq\bm G_n\mid X_s,X_t) \leq \mathbb P(\mathcal A_{w, \{s,t\}}\mid X_s,X_t) = \mathbb P(\mathcal A_{w, \{s,t\}}) = \frac{\lambda_{p_1,p_2,p_{12}}(\bm H)}{\lambda_{p_1,p_2,p_{12}}(\bm H_{w, \{s,t\} })}.
\end{align*}
Then recalling that $\mu = \lambda_{p_1, p_2, p_{12}}(\bm H)$, $X_w = I_w - \mu$, and using $|X_w| \leq I_w + \mu$, 
\begin{align*}
\mathbb E[|X_w| \mid X_s,X_t] \leq \frac{2  \lambda_{p_1, p_2, p_{12}}(\bm H) }{\lambda_{p_1,p_2,p_{12}}(\bm H_{w, \{s,t\} })}.
\end{align*} 
For every fixed submultiplex $\bm F\subseteq\bm H$, the number
of copies $\bm H_w$ satisfying $\bm H_{w, \{s,t\} }\cong\bm F$
is $\lesssim_{\bm H} n^{|V(\bm H)|-|V(\bm F)|}$.  Therefore, uniformly in $s,t$,
\begin{align*}
\sum_{w\in\overline N_{\mathcal G}(s,t)}
\mathbb E[|X_w|\mid X_s,X_t] \lesssim_{\bm H} n^{|V(\bm H)|} \lambda_{p_1, p_2, p_{12}}(\bm H)  \sum_{\substack{\bm F\subseteq\bm H\\|E(\bm F)|\geq1}}
\frac{1}{ n^{|V(\bm F)|} \lambda_{p_1,p_2,p_{12}}(\bm F)} \lesssim_{\bm H}
\frac{\mathbb E[X(\bm H,\bm G_n)]}{\Phi_{\bm H}} , 
\end{align*}
where the last step uses \eqref{eq:XHKnp}. Now, recall from \eqref{eq:Hvariance} that $
    \sigma^2 = \Var [ X( \bm{H}, \bm G_n)] \asymp_{\bm H} \frac{ (\E [X(\bm H, \bm G_n)])^2}{ \Phi_{\bm{H}}}$, under the parametrization \eqref{eq:exponent}. Combining the above with Theorem \ref{thm:W} gives the result in \eqref{eq:ZHGn}, as required. \hfill $\Box$

\section{Proof of Theorem \ref{thm:balancedH}}
\label{sec:balanceHpf}

We will prove Theorem \ref{thm:balancedH} using the Stein-Chen method for Poisson approximation based on dependency graphs.

\begin{theorem}[{\cite[Theorem~15]{poissonapproximation}}]  
Suppose $\{X_\alpha:\alpha\in A\}$ be collection of Bernoulli random variables and 
$X:=\sum_{\alpha\in A} X_\alpha$. Let $\mathcal G$ be a dependency graph for this collection.  Then, with $\lambda:=\mathbb E[X]$, the following holds: 
$$d_{\mathrm{TV}} \left(X,\operatorname{Pois}(\lambda) \right) \leq
\min\left\{1, \frac{1}{\lambda} \right\} \left( B_1+B_2 \right) , $$ 
where 
$$B_1 := \sum_{\alpha\in A} \sum_{\beta\in \overline{N}_{\mathcal{G}}(\alpha)}
\mathbb E[X_\alpha]\mathbb E[X_\beta] \text{ and } B_2 :=
\sum_{\alpha \in A}
\sum_{\substack{\beta\in \overline{N}_{\mathcal{G}}(\alpha)\\\beta\neq\alpha}}
\mathbb E[X_\alpha X_\beta].$$ 
(Here, $d_{\mathrm{TV}}$ denotes the total variation distance between two random variables.) 
\label{thm:poissondependencygraph} 
\end{theorem}

As in the proof of Theorem \ref{thm:ZHGn}, define  $I_s:=\bm{1}\{ \bm H_s \subseteq \bm G_n\}$, for $1 \leq s \leq M$. Recall, from \eqref{eq:XHKnp}, that $\mu := \mathbb E[I_s] = \lambda_{p_1,p_2,p_{12}}(\bm H)= n^{- |V(\bm H)|}, $  
where the last equality uses $\ell_{\bm H}(\bm\theta)=0$,  since $\bm H$ is $\bm\theta$-balanced. This implies that $\lambda := \mathbb E[X(\bm H,\bm G_n)] \rightarrow \frac{1}{|\mathrm{Aut}(\bm H)|}$. Construct a dependency graph $\mathcal G$ of the collection of random variables $\{I_s : 1 \leq s \leq M \}$, as in the proof of Theorem \ref{thm:ZHGn}. 

To apply Theorem \ref{thm:poissondependencygraph} first note that, for any $1 \leq s \leq M$, one has $|\overline{N}_{\mathcal{G}}(s)| \lesssim_{\bm H}
n^{|V(\bm H)|-2}$. Then, since $M \asymp_{\bm H} n^{|V(\bm H)|}$, 
\begin{align}
B_1= \sum_{s=1}^M \sum_{t \in \overline{N}_{\mathcal{G}}(s)} \mu^2 \lesssim_{\bm H} \frac{1}{n^2}.
\label{eq:B1}
\end{align}
It remains to estimate $B_2$. For this, fix two distinct adjacent copies $\bm H_s$ and $\bm H_t$ in the dependency graph and denote by $F_{\{s, t\}}$ the (one-sided) overlap multiplex as defined in the proof of Lemma~\ref{lm:meanvariance}. Then from \eqref{eq:GnHst}
\begin{equation}\label{eq:poissonst}
\mathbb P(I_s=I_t=1) \leq \frac{\mu^2}{\lambda_{p_1,p_2,p_{12}}
(\bm F_{\{s, t\}})}.
\end{equation}
We now consider the following two cases:  \\ 

\noindent
\emph{Case} $1$: $\bm F_{\{s, t\}}$ is a proper submultiplex of
$\bm H$. In this case, the number of ordered pairs $(s, t)$ satisfying $\bm F_{\{s, t\}}\cong\bm F$ is at most $\lesssim_{\bm H} n^{2 |V(\bm H)|-|V(\bm F)|}$.  Therefore, using \eqref{eq:poissonst}, the total
contribution to $B_2$ of all such pairs is at most
\begin{align}
\lesssim_{\bm H}
\sum_{\substack{ \bm F\subsetneq\bm H\\ |E(\bm F)|\geq1 }} n^{2 |V(\bm H)|-|V(\bm F)|}
\frac{\mu^2}{ \lambda_{p_1,p_2,p_{12}}(\bm F) }  = \sum_{\substack{
\bm F\subsetneq\bm H\\ |E(\bm F)|\geq1 }} n^{-\ell_{\bm F}(\bm\theta)} \ll 1, 
\label{eq:B2proper}
\end{align}
where the last step follows by noting that $\ell_{\bm F}(\bm\theta)>0$ for every proper submultiplex $\bm F\subsetneq\bm H$ containing an edge, since $\bm H$ is strictly $\bm\theta$-balanced. \\

\noindent  \emph{Case} $2$: $\bm F_{\{s, t\}}\cong\bm H$. In this case, the two copies have the same vertex set and the same underlying edge support, but since $s\neq t$, their layer assignments cannot be identical. This means that relative to one of the two copies, the joint occurrence requires an additional layer on at least one underlying edge. Consequently, for at least one underlying edge, the joint occurrence of the two copies incurs an additional factor of either $p_{12}/p_1$ or $p_{12}/p_2$ relative to the occurrence probability of one copy. This implies, $\mathbb P(I_s=I_t=1) \leq \mu n^{-\eta}$,  where $\eta = \min\left\{
\theta_{12}-\theta_1,\, \theta_{12}-\theta_2 \right\}>0$, by assumption. Moreover, for each $1\leq s\leq M$, there are only $\lesssim_{\bm H} 1 $ choices of $t$, having the same vertex set and underlying edge support as $\bm H_s$. Therefore, the total contribution to $B_2$ in this case is at most
\begin{align}
\lesssim_{\bm H} n^{|V(\bm H)|} \mu n^{-\eta} \lesssim_{\bm H} n^{-\eta} \ll 1.  
\label{eq:B2same}
\end{align}
Combining \eqref{eq:B1}, \eqref{eq:B2proper}, and \eqref{eq:B2same}, with Theorem~\ref{thm:poissondependencygraph}, the result in Theorem \ref{thm:balancedH} follows.  \hfill $\Box$

\section{Proof of Theorem \ref{thm:unbalancedH}}
\label{sec:unbalancedHpf}

Fix a point $\bm \theta= (\theta_1, \theta_2, \theta_{12}) \in \partial_{\bm H}$. Recall that,  $\mathcal{A}_{\bm H}(\bm \theta)$ denotes the collection of $\bm \theta$-extremal submultiplexes of $\bm H$. We begin by showing that the core $\overline{\bm{H}}$ of $\bm{H}$, which is the maximal element of $\mathcal{A}_{\bm H}(\bm \theta)$ under containment, is well-defined. Towards this, we need the notion of completion of a multiplex. Throughout, whenever we consider a submultiplex $\bm F$ of $\bm H$, we will implicitly assume that $|E(\bm F)| \geq 1$.

\begin{definition}
    Let $\bm F \subseteq \bm H$ be a submultiplex. Define the \emph{completion} $\hat{\bm F} = (V(\bm F), \hat{F}^{(1)}, \hat{F}^{(2)}), $ of $\bm F = (V(\bm F), F^{(1)}, F^{(2)})$ as the multiplex with,\footnote{For a finite set $U$, let $\binom{U}{2} =
\{\{u,v\}: u,v\in U, u\neq v \} ,$
denote the collection of all unordered pairs of distinct elements
of $U$. }
\begin{equation*}
    E(\hat{F}^{(i)}) =   E(F^{(i)}) \cup \left(E(H^{(1)}) \cap E(H^{(2)}) \cap \binom{V\left(\bm F\right)}{2} \right) , 
\end{equation*} 
for $i \in \{1, 2\}$. If $\bm F=\hat{\bm F}$, then we call the submultiplex \emph{complete}.
\end{definition}


With this definition we can show that the core of $\bm H$ is well-defined.

\begin{lemma}
    \label{lem:extreme}
    If $\bm \cQ, \bm \cR \in \mathcal{A}_{\bm H}(\bm \theta)$, then $\hat{\bm \cQ}\cup \hat{\bm \cR} \in \mathcal{A}_{\bm H}(\bm \theta)$. Consequently, the maximum element in $\mathcal{A}_{\bm H}(\bm \theta)$ and the core $\overline{\bm H}$ is well-defined.
\end{lemma}

\begin{proof} We begin with the following fact (which is proved later in Section \ref{sec:doublepresentpf}):
    \begin{align}\label{eq:doublepresent}
        \bm F \in  \mathcal{A}_{\bm H}(\bm \theta) \Longrightarrow \bm F \subseteq \hat{\bm F} \in \mathcal{A}_{\bm H}(\bm \theta) .  
    \end{align}
    Hence, without loss of generality, we may assume that $\bm \cQ = \hat{\bm \cQ}$ and $\bm \cR = \hat{\bm \cR}$. Now, recall the definition of the function $\ell_{\bm F}(\bm \theta)$ from \eqref{eq:VF}. The key observation is that the function  $\ell_{\bm F}(\bm \theta)$ is {\it modular}, that is,
    \begin{align}\label{eq:modular}
        \ell_{\cQ\cup \cR}(\bm{\theta})=\ell_{\cQ}(\bm{\theta})+\ell_{\cR}(\bm{\theta})- \ell_{\cQ\cap \cR}(\bm{\theta}) ,  
    \end{align} 
     for $\bm \cQ, \bm \cR \in \mathcal{A}_{\bm H}(\bm \theta) $ which are both complete. We prove this fact in Section \ref{sec:modularpf} below. This implies, since $\ell_{\cQ}(\bm{\theta}) = \ell_{\cR}(\bm{\theta}) = 0$ (recall $\cQ,  \cR$ are $\bm \theta$-extremal),
    \begin{align}\label{eq:F12}
        \ell_{\cQ\cup \cR}(\bm{\theta})= - \ell_{\cQ\cap \cR}(\bm{\theta}) .
    \end{align}
    Also, since $\bm \theta  \in \partial_{\bm H}$, we have $\ell_{\bm{F}}(\bm{\theta}) \geq 0$ for all submultiplexes $\bm F \subseteq \bm H$. This means, $\ell_{\cQ\cap \cR}(\bm{\theta}) \geq 0$ and $\ell_{\cQ\cup \cR}(\bm{\theta}) \geq 0$, since $\cQ, \cR$ are both submultiplexes of $\bm{H}$. Hence, \eqref{eq:F12} can hold only if $\ell_{\cQ\cup \cR}(\bm{\theta}) = 0$, which means $\cQ\cup \cR$ is a $\bm \theta$-extremal submultiplex of $\bm H$.  
\end{proof}

\subsection{Proof of \eqref{eq:doublepresent}} \label{sec:doublepresentpf}

We begin with the following lemma:

\begin{lemma}
    \label{monotonicity}
    Fix a multiplex $\bm F_1 = (V(\bm F_1), F_1^{(1)}, F_1^{(2)})$. Define $\bm F_2=\bm F_1\cup \bm{\mathcal{O}}_{(u, v)}$, where $\bm{\mathcal{O}}_{(u, v)}$ is a multiplex with vertex set $\{u, v\} \subseteq V(\bm F_1)$ with the edge $(u, v) $ present in both layers.  
    Then $\ell_{\bm F_1}(\bm \theta)\geq \ell_{\bm F_2}(\bm \theta)$, for $\bm \theta \in \Theta$.
\end{lemma}

\begin{proof}  
    If $(u,v) \in E(F_1^{(1)}) \cap E(F_1^{(2)})$, then $\bm F_1=\bm F_2$ and $\ell_{\bm F_1}(\bm \theta)= \ell_{\bm F_2}(\bm \theta)$. Otherwise, there are the following three possibilities:

    \begin{itemize}
        \item[(1)] $(u,v)\notin E(F_1^{(1)})\cup E(F_1^{(2)})$: Then, as $\theta_{12} > 0 $,
            \begin{align*}
                \ell_{\bm F_2}(\bm \theta) =|V(\bm{F_1})|- \theta_1 |E(F_1^{(1)})\setminus  E(F_1^{(2)})| & -  \theta_2 |E(F_1^{(2)})\setminus E(F_1^{(1)})|                                  \\
                                                                                                          & - \theta_{12} (|E(F_1^{(1)})\cap E(F_1^{(2)})|+1) < \ell_{\bm{F_1}}(\bm \theta) .
            \end{align*}

        \item[(2)] $(u,v)\in E(F_1^{(1)})\setminus E(F_1^{(2)})$: Then, as $\theta_{12}\geq\theta_1$,
            \begin{align*}
                \ell_{\bm F_2}(\bm \theta) =|V(\bm{F_1})|- \theta_1 (|E(F_1^{(1)}) \setminus  E(F_1^{(2)})| & -1)   - \theta_2 |E(F_1^{(2)})\setminus E(F_1^{(1)})|                             \\
                                                                                                            & - \theta_{12} (|E(F_1^{(1)})\cap E(F_1^{(2)})|+1)\leq\ell_{\bm F_1}(\bm \theta) .
            \end{align*}

        \item[(3)] $(u,v)\in E(F_1^{(2)})\setminus E(F_1^{(1)})$: Then,  as $\theta_{12}\geq\theta_2$,
            \begin{align*}
                \ell_{\bm F_2}(\bm \theta) = |V(\bm{F_1})|- \theta_1 |E(F_1^{(1)})\setminus E(F_1^{(2)})| & - \theta_2 (|E(F_1^{(2)})\setminus E(F_1^{(1)})|-1)                                 \\
                                                                                                          & - \theta_{12} (|E(F_1^{(1)})\cap E(F_1^{(2)})|+1) \leq \ell_{\bm F_1}(\bm \theta) .
            \end{align*}

    \end{itemize}
    This completes the proof of Lemma \ref{monotonicity}.
\end{proof}

Now, suppose $\bm F \in \mathcal{A}_{\bm H}(\bm \theta)$ and $\hat{\bm F}$ is the completion of $\bm F$. This means, by Lemma \ref{monotonicity}, $0 \leq  \ell_{\hat{\bm F}}(\bm \theta)\leq \ell_{\bm F}(\bm \theta) = 0$. Hence, $\ell_{\hat{\bm F}}(\bm \theta) = 0$, that is, $\hat{\bm F} \in \mathcal{A}_{\bm H}(\bm \theta)$.

\subsection{Proof of \eqref{eq:modular}} \label{sec:modularpf}

Let $\bm Q = (V(\bm Q), Q^{(1)}, Q^{(2)}) \text { and } \bm R = (V(\bm R), R^{(1)}, R^{(2)})$ be two complete submultiplexes in $\mathcal{A}_{\bm H}(\bm \theta)$. Observe that both sides of \eqref{eq:modular} are linear functions of $\theta_1, \theta_2, \theta_{12}$. To begin with, observe that $
    |V(\bm \cQ \cup \cR)| = |V(\bm \cQ)| + |V(\cR)| -  |V(\bm \cQ \cap \cR)|$. Hence, the constant terms on both sides of \eqref{eq:modular} are equal. Next, consider the coefficients of $\theta_1$ on both sides of \eqref{eq:modular}. To this end, note that
\begin{align*}
      & \left(E(Q^{(1)}) \cup E(R^{(1)})\right) \setminus \left(E(Q^{(2)}) \cup E(R^{(2)})\right)                                              \\
    & = \left(E(Q^{(1)}) \cup E(R^{(1)})\right) \cap \left( E(Q^{(2)})^{c} \cap E(R^{(2)})^{c} \right)                                         \\
    & = \left(E(Q^{(1)}) \cap E(Q^{(2)})^{c} \cap E(R^{(2)})^{c}\right) \cup \left(E(R^{(1)}) \cap E(Q^{(2)})^{c} \cap E(R^{(2)})^{c}\right) .
\end{align*}
Observe, if $(u,v) \in E(Q^{(1)}) \cap E(Q^{(2)})^{c}$, then $(u,v)$ is an edge only in layer 1 of $\bm \cQ$. Therefore, by completeness of $\bm \cQ$, $(u,v)$ can only be an edge in one layer of $\bm H$. This means $(u,v) \in E(R^{(2)})^c$. Thus, $E(Q^{(1)}) \cap E(Q^{(2)})^{c} \cap E(R^{(2)})^{c} = E(Q^{(1)}) \cap E(Q^{(2)})^{c}$. Similarly, $E(R^{(1)}) \cap E(Q^{(2)})^{c} \cap E(R^{(2)})^{c} = E(R^{(1)}) \cap E(R^{(2)})^{c}$. Combining the above gives,
\begin{equation*}
    \left(E(Q^{(1)}) \cup E(R^{(1)})\right) \setminus \left(E(Q^{(2)}) \cup E(R^{(2)})\right) =  \left(E(Q^{(1)}) \setminus E(Q^{(2)})\right) \cup \left(E(R^{(1)}) \setminus E(R^{(2)})\right) .
\end{equation*}
Hence,
\begin{align}\label{eq:theta1}
      & \left|\left(E(Q^{(1)}) \cup E(R^{(1)})\right) \setminus \left(E(Q^{(2)}) \cup E(R^{(2)})\right)\right| \nonumber                                                                                            \\
    &  =  \left|E(Q^{(1)}) \setminus E(Q^{(2)})\right|+ \left|E(R^{(1)}) \setminus E(R^{(2)})\right| -  \left|\left(E(Q^{(1)}) \setminus E(Q^{(2)})\right) \cap \left(E(R^{(1)}) \setminus E(R^{(2)})\right)\right| .  
\end{align}
Note that
\begin{align*}
    (E(Q^{(1)}) \setminus E(Q^{(2)})) \cap (E(R^{(1)}) \setminus E(R^{(2)})) & =  (E(Q^{(1)})\cap E(R^{(1)}))\setminus (E(Q^{(2)})\cup E(R^{(2)})) \nonumber \\
                                                                             & \subseteq (E(Q^{(1)})\cap E(R^{(1)}))\setminus (E(Q^{(2)})\cap E(R^{(2)})) .
\end{align*}
We will show that the inclusion in the other direction also holds. For this, suppose $$(u,v) \in (E(Q^{(1)})\cap E(R^{(1)}))\setminus (E(Q^{(2)})\cap E(R^{(2)})).$$ Then $(u,v)\in E(Q^{(1)})\cap E(R^{(1)})$ and $(u,v)\in E(Q^{(2)})^c \cup E(R^{(2)})^c$. Hence, there are two possibilities:
\begin{itemize}
    \item $(u,v)\notin E(Q^{(2)})$. This means $(u,v) \in E(Q^{(1)}) \setminus E(Q^{(2)})$. Also, by completeness of $\bm \cQ$, $(u,v)\notin E(H^{(2)})$, and, as a result, $(u,v) \notin E(R^{(2)})$. This shows, $(u,v) \in E(R^{(1)}) \setminus E(R^{(2)})$. Hence, $(u,v) \in  (E(Q^{(1)}) \setminus E(Q^{(2)})) \cap (E(R^{(1)}) \setminus E(R^{(2)}))$.

    \item $(u,v)\notin E(R^{(2)})$. In this case, similarly, it can be shown that $(u,v) \in  (E(Q^{(1)}) \setminus E(Q^{(2)})) \cap (E(R^{(1)}) \setminus E(R^{(2)}))$.
\end{itemize}
Hence, combining the above shows,
\begin{align*}
    (E(Q^{(1)}) \setminus E(Q^{(2)})) \cap (E(R^{(1)}) \setminus E(R^{(2)})) & =  (E(Q^{(1)})\cap E(R^{(1)}))\setminus (E(Q^{(2)})\cup E(R^{(2)})) \nonumber \\
                                                                             & = (E(Q^{(1)})\cap E(R^{(1)}))\setminus (E(Q^{(2)})\cap E(R^{(2)})) .
\end{align*}
Therefore, from \eqref{eq:theta1},
\begin{align*}
      & \left|\left(E(Q^{(1)}) \cup E(R^{(1)})\right) \setminus \left(E(Q^{(2)}) \cup E(R^{(2)})\right)\right|                                                                                                 \\
    = & \left|E(Q^{(1)}) \setminus E(Q^{(2)})\right|+ \left|E(R^{(1)}) \setminus E(R^{(2)})\right| -  \left|\left(E(Q^{(1)}) \cap E(R^{(1)})\right) \setminus \left(E(Q^{(2)}) \cap E(R^{(2)})\right)\right| .
\end{align*}
This shows that the coefficients of $\theta_1$ on both sides of \eqref{eq:modular} are equal. Similarly, it can be shown that the coefficients of $\theta_2$ on both sides of \eqref{eq:modular} are equal. Finally, consider the coefficients $\theta_{12}$ on both sides of \eqref{eq:modular}. Towards this,
\begin{align}\label{eq:QR}
     & (E(Q^{(1)}) \cup E(R^{(1)})) \cap (E(Q^{(2)}) \cup E(R^{(2)}))                                                                                          \\
     & = (E(Q^{(1)}) \cap E(Q^{(2)}) ) \cup (E(R^{(1)}) \cap E(R^{(2)}) )  \cup (E(Q^{(1)}) \cap E(R^{(2)}) )   \cup (E(Q^{(2)}) \cap E(R^{(1)}) ) . \nonumber
\end{align}
Note that
\begin{align}\label{eq:QR12}
    (E(Q^{(1)}) \cap E(Q^{(2)}) ) \cup (E(Q^{(1)}) \cap E(R^{(2)}) )   = E(Q^{(1)}) \cap E(Q^{(2)}) .
\end{align}
Clearly, the RHS is contained in the LHS. Now, suppose $(u, v) \in (E(Q^{(1)}) \cap E(R^{(2)})) \setminus ( E(Q^{(1)}) \cap E(Q^{(2)}) )$. This means $(u, v) \in (E(Q^{(1)}) \cap E(R^{(2)}))$ and $ (u, v) \notin E(Q^{(2)})$. Then by completeness of $\bm \cQ$, $ (u, v) \notin E(H^{(2)})$, and as a result $(u, v) \notin E(R^{(2)})$, which is a contradiction. Hence, $(E(Q^{(1)}) \cap E(R^{(2)})) \setminus ( E(Q^{(1)}) \cap E(Q^{(2)}) ) = \emptyset$, and \eqref{eq:QR12} holds. Similarly,
\begin{align}\label{eq:RQ12}
    (E(R^{(1)}) \cap E(R^{(2)}) ) \cup (E(Q^{(2)}) \cap E(R^{(1)}) )   = E(R^{(1)}) \cap E(R^{(2)}) .
\end{align}
Combining \eqref{eq:QR}, \eqref{eq:QR12}, and \eqref{eq:RQ12} gives,
$$(E(Q^{(1)}) \cup E(R^{(1)})) \cap (E(Q^{(2)}) \cup E(R^{(2)})) = (E(Q^{(1)}) \cap E(Q^{(2)}) ) \cup (E(R^{(1)}) \cap E(R^{(2)}) )  . $$
Hence,
\begin{align*}
    |(E(Q^{(1)}) \cup E(R^{(1)})) \cap (E(Q^{(2)}) \cup E(R^{(2)}))| & = | E(Q^{(1)}) \cap E(Q^{(2)})|  + | E(R^{(1)}) \cap E(R^{(2)}) | \nonumber \\
                                                                     & - | E(Q^{(1)}) \cap E(Q^{(2)}) \cap E(R^{(1)}) \cap E(R^{(2)})| . \nonumber
\end{align*}
This shows that the coefficients of $\theta_{12}$ on both sides of \eqref{eq:modular} are equal.  \hfill $\Box$

\subsection{ Completing the Proof of Theorem \ref{thm:unbalancedH} }

Note that since $\overline{\bm H }$ is the core, $\ell_{\overline{\bm H }}(\bm \theta) = 0$. 
We begin by stating some consequences of the definition of the core.  

\begin{itemize} 

\item By \eqref{eq:doublepresent}, the completion
$\hat{\overline{\bm H}}$ is also $\bm\theta$-extremal. Since ${\overline{\bm H}}$ is the maximum $\bm\theta$-extremal submultiplex of $\bm H$ and
${\overline{\bm H}}\subseteq\hat{\overline{\bm H}}$, it follows that
$\hat{\overline{\bm H}}={\overline{\bm H}}$. Thus, the core is
complete.

\item Moreover, ${\overline{\bm H}}$ is the unique copy of its isomorphism
type contained in $\bm H$. Indeed, suppose that
$\bm F\subseteq\bm H$ and
$\bm F\cong{\overline{\bm H}}$. Then
$\ell_{\bm F}(\bm\theta)=0$, so $\bm F$ is
$\bm\theta$-extremal. Since ${\overline{\bm H}}$ is the maximum
extremal submultiplex, $\bm F\subseteq{\overline{\bm H}}$.
The multiplexes $\bm F$ and ${\overline{\bm H}}$ are isomorphic,
and therefore have the same number of vertices and the same edge
counts in both layers. Hence, the containment
$\bm F\subseteq{\overline{\bm H}}$ forces
$\bm F={\overline{\bm H}}$. Consequently, every distinct image copy of $\bm H$ contains
exactly one distinct image copy of ${\overline{\bm H}}$.

\item If ${\overline{\bm H}}=\bm H$, then
$\ext_n({\overline{\bm H}},\bm H)=1$ almost surely, and the quantity
in \eqref{eq:unbalancedH} is identically zero. We may therefore
assume that ${\overline{\bm H}}\subsetneq\bm H$.
\end{itemize}
To proceed, fix a copy $\overline{\bm H}'$ of ${\overline{\bm H}}$ in the
complete multiplex $\bm K_n$. Let
$\mathcal N_n(\overline{\bm H}')$ denote the collection of all distinct copies $\bm{H}'$ of $\bm H$ in $\bm K_n$ satisfying
$\overline{\bm H}'\subseteq\bm{H}'$. For
$\bm{H}'\in\mathcal N_n(\overline{\bm H}')$, define
$$Y_{\bm{H}'} := \bm 1 \{\bm{H}'\subseteq\bm G_n\cup\overline{\bm H}' \} \text{ and  } Z_{\overline{\bm H}'} := \sum_{\bm{H}'\in\mathcal N_n(\overline{\bm H}')}Y_{\bm{H}'}.$$ Then $Z_{\overline{\bm H}'}$ has the same distribution as
$\ext_n({\overline{\bm H}},\bm H)$. By exchangeability, this
distribution does not depend on the particular copy $\overline{\bm H}'$. Also, let
$\mathcal C_n({\overline{\bm H}})$ be the collection of distinct
copies of the core in $\bm K_n$. Note that on the event $\{
\overline{\bm H}'\subseteq\bm G_n
\}$, one has
$\bm G_n\cup\overline{\bm H}'=\bm G_n$. Hence,
$\bm 1\{
\overline{\bm H}'\subseteq\bm G_n
\}Z_{\overline{\bm H}'}$ counts the copies of $\bm H$ occurring
in $\bm G_n$ that contain the particular core copy
$\overline{\bm H}'$. Since every copy of $\bm H$ contains exactly one
copy of ${\overline{\bm H}}$, 
\begin{equation}\label{eq:extensionH}
X(\bm H,\bm G_n) = \sum_{\overline{\bm H}'\in\mathcal C_n({\overline{\bm H}})}
\bm 1\{ \overline{\bm H}'\subseteq\bm G_n \}Z_{\overline{\bm H}'} \text{ and } X({\overline{\bm H}},\bm G_n) = \sum_{\overline{\bm H}'\in\mathcal C_n({\overline{\bm H}})} \bm 1\{ \overline{\bm H}'\subseteq\bm G_n \}.
\end{equation}   
Note that the random variables $\bm 1\{\overline{\bm H}'\subseteq\bm G_n
\}$ and $Z_{\overline{\bm H}'}$ are independent. This is because, $\bm 1\{
\overline{\bm H}'\subseteq\bm G_n\}$ depends only on the edge
vectors corresponding to the underlying edges of $\overline{\bm H}'$. Also, since ${\overline{\bm H}}$ is complete, an extension does not impose
an additional layer requirement on any of these underlying
edges. After the edges of $\overline{\bm H}'$ have been planted,
$Z_{\overline{\bm H}'}$ therefore depends only on edge vectors indexed by
underlying pairs outside $E(\overline{\bm H}')$. Hence, using \eqref{eq:extensionH}, the triangle inequality, independence, and the Cauchy--Schwarz inequality, we
obtain
\begin{align}
\mathbb E\left[
\left|\frac{X(\bm H,\bm G_n)}{ \E[\ext_n({\overline{\bm H}},\bm H)]} - X({\overline{\bm H}},\bm G_n) \right| \right] &\leq \sum_{\overline{\bm H}'\in\mathcal C_n({\overline{\bm H}})}
\mathbb E\left[\bm 1\{\overline{\bm H}'\subseteq\bm G_n\} \left| \frac{Z_{\overline{\bm H}'}}{ \E[Z_{\overline{\bm H}'}] }-1 \right| \right] \nonumber  \\
&=
\mathbb E[X({\overline{\bm H}},\bm G_n)] \mathbb E\left[ \left| \frac{Z_{\overline{\bm H}'}}{ \E[Z_{\overline{\bm H}'}] }-1 \right| \right]
\nonumber\\
&\leq
\mathbb E[X({\overline{\bm H}},\bm G_n)]
\sqrt{ \frac{ \mathrm{Var}[Z_{\overline{\bm H}'}]}{ \E[Z_{\overline{\bm H}'}]^2  } } . 
\label{eq:extension-L1-bound}
\end{align}
Note that $\mathbb E[X({\overline{\bm H}},\bm G_n)]
\asymp_{\bm H} n^{|V(\overline{\bm H})|} \lambda_{p_1, p_2, p_{12}}({\overline{\bm H}}) = n^{\ell_{{\overline{\bm H}}}(\bm\theta)} = 1$. Hence, to show \eqref{eq:unbalancedH}  it suffices to prove that 
\begin{align}\label{eq:extensionsecondmoment}
\mathrm{Var}[Z_{\overline{\bm H}'}] \ll \E[Z_{\overline{\bm H}'}]^2  .  
\end{align} 
%
%
%
%

\subsubsection{Proof of \eqref{eq:extensionsecondmoment} }

Since the core is complete, an extension of $\overline{\bm H}'$ does not
require an additional layer on any underlying edge already used
by $\overline{\bm H}'$. Hence, for every
$\bm{H}'\in\mathcal N_n(\overline{\bm H}')$, $\mathbb P(Y_{\bm{H}'}=1) = \frac{\lambda_{p_1, p_2, p_{12}}(\bm H)} {\lambda_{p_1, p_2, p_{12}}({\overline{\bm H}})}$. Furthermore, $|\mathcal N_n(\overline{\bm H}')|\asymp_{\bm H} n^{|V(\bm H)|-|V(\overline{\bm H})}$.
Hence, it follows that
\begin{equation}\label{eq:extension-mean}
\E[Z_{\overline{\bm H}'}] \asymp_{\bm H}
n^{|V(\bm H)| - |V(\overline{\bm H})}
\frac{\lambda_{p_1, p_2, p_{12}}(\bm H)}
{\lambda_{p_1, p_2, p_{12}}({\overline{\bm H}})}
=
n^{\ell_{\bm H}(\bm\theta)
-\ell_{{\overline{\bm H}}}(\bm\theta)}
=
n^{\ell_{\bm H}(\bm\theta)}.
\end{equation}
We now compute the variance of $Z_{\overline{\bm H}'}$. For this, for an extension $\bm{H}'$, define its residual underlying edge
support by $S(\bm{H}') := E(\bm{H}')\setminus E(\overline{\bm H}')$.  For two extensions $\bm{H}',\bm H''$, if $S(\bm{H}')\cap S(\bm H'')=\varnothing$, then
$Y_{\bm{H}'}$ and $Y_{\bm H''}$ depend on disjoint collections
of independent edge vectors and are therefore independent.
For two extensions $\bm{H}',\bm H''$ satisfying
$S(\bm{H}')\cap S(\bm H'')\neq\varnothing$, define the (one-sided) 
overlap multiplex $\bm F_{\bm{H}', \bm H''}$ by $V(\bm F_{\bm{H}', \bm H''}):=
V(\bm{H}')\cap V(\bm H'')$, and, for $i\in\{1,2\}$,
\[
E(F_{\bm{H}', \bm H''}^{(i)})
:=
E((\overline{H}')^{(i)})
\cup
\left[
\left(
E((H')^{(i)})\setminus E((\overline{H}')^{(i)})
\right)
\cap E(\bm H'')
\right].
\]
An argument similar to the one used in \eqref{eq:GnHst} gives
\begin{equation}\label{eq:extension-joint-bound}
\mathbb P(Y_{\bm{H}'}=Y_{\bm H''}=1)
\leq
\frac{
\left(
\lambda_{p_1, p_2, p_{12}}(\bm H)/\lambda_{p_1, p_2, p_{12}}({\overline{\bm H}})
\right)^2
}{
\lambda_{p_1, p_2, p_{12}}(\bm F_{\bm{H}', \bm H''})
/
\lambda_{p_1, p_2, p_{12}}({\overline{\bm H}})
}.
\end{equation}
For every fixed isomorphism type $\bm F$ arising in this way,
the number of ordered pairs $(\bm{H}',\bm H'')$ satisfying
$\bm F_{\bm{H}', \bm H''}\cong\bm F$ is at most $\lesssim_{\bm H} n^{2(|V(\bm H)| - |V(\overline{\bm H})| )-(|V(\bm F)| -  |V(\overline{\bm H})|)}$.  
Let $\mathcal  F$ be the finite collection of isomorphism
types of submultiplexes of $\bm H$ that properly contain a
copy of ${\overline{\bm H}}$ and can arise as such an overlap. Then, using $\mathrm{Cov}[Y_{\bm{H}'},Y_{\bm H''}] \leq \mathbb P(Y_{\bm{H}'}=Y_{\bm H''}=1)$, \eqref{eq:extension-mean}, and \eqref{eq:extension-joint-bound} gives, 
\begin{align*}
\mathrm{Var}[Z_{\overline{\bm H}'}] & \lesssim_{\bm H}\, n^{2( |V(\bm H)| - |V(\overline{\bm H})|)}
\left(\frac{\lambda_{p_1, p_2, p_{12}}(\bm H)} {\lambda_{p_1, p_2, p_{12}}({\overline{\bm H}})} \right)^2  \sum_{\bm F\in\mathcal  F} \frac{ \lambda_{p_1, p_2, p_{12}}({\overline{\bm H}}) }{ n^{|V(\bm F)|- |V(\overline{\bm H})| } \lambda_{p_1, p_2, p_{12}}(\bm F) }  \nonumber \\  
& \lesssim_{\bm H} \E[Z_{\overline{\bm H}'}]^2  \sum_{\bm F\in\mathcal  F} \frac{1}{n^{\ell_{\bm F}(\bm\theta)}} \ll \E[Z_{\overline{\bm H}'}]^2 ,  
\end{align*}
since $\ell_{\bm F}(\bm\theta) > 0$, for all $\bm F\in\mathcal  F$, by the maximality of the core. This proves \eqref{eq:extensionsecondmoment}.  \hfill $\Box$

\color{black}

\section{ Discussion and Future Directions }
\label{sec:summary}

In this paper we have initiated the study of thresholds for submultiplexes and their fluctuations in correlated multiplex networks. Although we have considered  only two layers for simplicity, most of the results can be naturally extended to more than two layers. Since the key ideas remain the same, we have chosen to present our results for the two-layer case to keep the exposition accessible. Multiplex networks can also be viewed within the broader framework of colored or weighted graphs (see \cite{multiplexnetworks,zucal2024probabilitygraphons,abraham2025probability} and the references therein).  This suggests potential connections to broader topics in combinatorics and opens up possibilities for many natural generalizations that merit further investigation.

The following are a few more possible directions for further developing the study of submultiplexes in correlated networks.

\begin{itemize}

    \item In this paper we have derived the threshold for the emergence of a `small' (fixed size) multiplex in the correlated Erd\H{o}s-R\'enyi multiplex model. However, sometimes one is also interested in `large' submultiplexes that grow with the size of the network. In the context of random graphs, there are many celebrated results identifying thresholds for various types of large subgraphs, particularly for spanning structures such as Hamiltonian cycles and perfect matchings (see the survey \cite{perkins2025searching}). Determining the corresponding thresholds for analogous global structures in the multiplex setting could be an interesting future direction. In this context, it might also be worth exploring the implications of the celebrated Kahn-Kalai conjecture \cite{kahn2007thresholds} (and the recent breakthrough \cite{conjectureproof}) in the multiplex/colored graph setting.

    \item Another direction would be to study the geometric properties of the satisfiability region. For example, given a submultiplex $\bm{H}$, it would be interesting to investigate how the number of edges or vertices (extreme points) of $\partial_{\bm H}$ depends on the structure of $\bm{H}$.

\end{itemize}

\small

\subsection*{Acknowledgements} This work was carried out in part during the Workshop on the Mathematical Foundations of Network Models and Their Applications, held as part of the BIRS-CMI pilot program organized by Louigi Addario-Berry, Siva Athreya, Shankar Bhamidi, Serte Donderwinkel, and Soumik Pal. The workshop took place at the Chennai Mathematical Institute (CMI) from December 15-20, 2024, and was preceded by a Research School held from December 9-13, 2024. The events were supported by Banff International Research Station (BIRS), Chennai Mathematical Institute (CMI), and the National Board for Higher Mathematics (NBHM). The authors thank the organizers for their inspirational leadership and hospitality. The authors also thank Siva Athreya for many useful discussions. The authors are also indebted to an anonymous referee for their careful comments, which improved the quality and the presentation of the paper. BBB was supported by NSF CAREER grant DMS 2046393 and a Sloan Research Fellowship. LE was supported by DGAPA-PAPIIT project IN108525. KD was supported by the Prime Minister's Research Fellowship PM/PMRF-22-18686.03 (PMRF ID: 0202965) from the Ministry of Education of India.

\bibliographystyle{abbrv}
\bibliography{bibliography,references}
\end{document}

%% file: plots/GraphMultiplex.tex
\begin{tikzpicture}[scale=1.1,
    every node/.style={circle, draw=black, fill=green!80, inner sep=2pt, font=\small}
]
\node[draw=none, fill=none, text width=1cm] at (3.5, 3.6) {$\bm{G}_n$};

\node (1)  at (0.5,2.8) {1};
\node (2)  at (0.5,1.3) {2};
\node (3)  at (-1.5,-0.2) {3};
\node (4)  at (-1.2,2.7) {4};
\node (5)  at (3.8,0.7) {5};
\node (6)  at (1.5,0) {6};
\node (7)  at (1.3,-1.0) {7};
\node (8)  at (4.0,1.9) {8};
\node (9)  at (0.0,0.2) {9};
\node (10) at (1.6,3.0) {10};
\node (11) at (3.0,2.3) {11};
\node (12) at (-1.3,0.8) {12};
\node (13) at (-1.4,1.9) {13};
\node (14) at (1.5,1.3) {14};
\node (15) at (-0.8,3.6) {15};

\draw[red, thick] (15) to[bend left=20] (4);
\draw[red, thick] (4) to[bend left=20] (13);
\draw[red, thick] (13) -- (12);
\draw[red, thick] (1) to[bend left=20] (10);
\draw[red, thick] (10) -- (11);
\draw[red, thick] (11) to[bend left=20] (8);
\draw[red, thick] (8) -- (5);
\draw[red, thick] (1) -- (2);
\draw[red, thick] (1) -- (15);
\draw[red, thick] (1) -- (13);
\draw[red, thick] (4) -- (9);
\draw[red, thick] (2) to[bend left=20] (9);
\draw[red, thick] (9) to[bend left=20] (6);
\draw[red, thick] (9) to[bend left=20] (3);
\draw[red, thick] (6) -- (7);
\draw[red, thick] (9) -- (14);

\draw[blue, thick] (9) -- (7);
\draw[blue, thick] (9) to[bend right=20] (3);
\draw[blue, thick] (15) to[bend right=20] (4);
\draw[blue, thick] (4) to[bend right=20] (13);
\draw[blue, thick] (1) to[bend right=20] (10);
\draw[blue, thick] (1) to[bend right=15] (11);
\draw[blue, thick] (11) to[bend right=20] (8);
\draw[blue, thick] (2) to[bend right=20] (9);
\draw[blue, thick] (9) to[bend right=20] (6);
\draw[blue, thick] (9) -- (12);
\draw[blue, thick] (12) -- (3);
\draw[blue, thick] (11) -- (6);
\draw[blue, thick] (14) -- (5);
\draw[blue, thick] (14) -- (2);
\draw[blue, thick] (14) -- (10);

\end{tikzpicture}

%% file: plots/Multiplex_1.tex
\begin{tikzpicture}[scale = 0.425]
  \begin{axis}[
    width=12cm,
    height=12cm,
    xmin=-0.05, xmax=1.10,
    ymin=-0.05, ymax=2.20,
    axis background/.style={fill=black!10},
    xlabel={\huge{$\theta$}},
    ylabel={\huge{$\theta_{12}$}},
  ]
    \addplot[
      fill=blue!40,
      draw=none,
      area legend
    ] coordinates {
      (0.000000, 0.000000)
      (0.000000, 2.000000)
      (0.500000, 2.000000)
      (1.000000, 1.000000)
    } -- cycle;
    \addplot[
      fill=white,
      draw=none
    ] coordinates {
      (-0.05, -0.05)
      (0, -0.05)
      (0, 2.20)
      (-0.05, 2.20)
    } -- cycle;
    \addplot[
      fill=white,
      draw=none
    ] coordinates {
      (0, 0)
      (0, -0.05)
      (1.10, -0.05)
      (1.10, 1.10)
    } -- cycle;
    \addplot[
      red,
      thick,
      domain=0.000000:0.500000
    ] {0.000000*(x-0.000000)+2.000000};
    \addplot[
      red,
      thick,
      domain=0.500000:1.000000
    ] {-2.000000*(x-0.500000)+2.000000};
    \addplot[
      black,
      dashed,
      domain=0:1.0
    ] {x};
    \addplot[
      black,
      dashed
    ] coordinates {
      (0, 0)
      (0, 2.0)
    };
    \addplot[
      only marks,
      mark=*,
      mark size=3pt,
      red
    ] coordinates {
    (1.000000, 1.000000)
      
      (0.500000, 2.000000)
    };
    \addplot[
      only marks,
      mark=*,
      mark size=3pt,
      black
    ] coordinates {
      (0.000000, 0.000000)
      (0.000000, 2.000000) 
      
    };
    
    \node at (axis cs:0.35,1.2) [font=\huge] {many copies of $\bm{\mathcal{R}}$};
    
    \node at (axis cs:0.85,1.9) [font=\huge] {no copies of $\bm{\mathcal{R}}$};
    
    \node at (axis cs:0.12,0.105) [anchor=north east, xshift=-3pt, yshift=-3pt] {\LARGE{$O$}};
    \node at (axis cs:0.090000,1.97500) [anchor=south east, xshift=-3pt, yshift=3pt] {\LARGE{$Q_3$}};
    \node at (axis cs:0.500000,1.95000) [anchor=south west, xshift=3pt, yshift=3pt] {\LARGE{$Q_2$}};
    \node at (axis cs:1.000000,1.0000) [anchor=west, xshift=3pt] {\LARGE{$Q_1$}};
    
  \end{axis}
\end{tikzpicture}

%% file: plots/Multiplex_2.tex
\begin{tikzpicture}[scale = 0.425]
  \begin{axis}[
    width=12cm,
    height=12cm,
    xmin=-0.05, xmax=1.4,
    ymin=-0.05, ymax=2.20,
    axis background/.style={fill=black!10},
    xlabel={\huge{$\theta$}},
    ylabel={\huge{$\theta_{12}$}},
  ]
    \addplot[
      fill=blue!40,
      draw=none,
      area legend
    ] coordinates {
      (0.000000, 0.000000)
      (0.000000, 2.000000)
      (0.500000, 2.000000)
      (0.666667, 1.666667)
      (0.833333, 0.833333)
    } -- cycle;
    \addplot[
      fill=white,
      draw=none
    ] coordinates {
      (-0.05, -0.05)
      (0, -0.05)
      (0, 2.20)
      (-0.05, 2.20)
    } -- cycle;
    \addplot[
      fill=white,
      draw=none
    ] coordinates {
      (0, 0)
      (0, -0.05)
      (2.2, -0.05)
      (2.2, 2.2)
    } -- cycle;
    \addplot[
      red,
      thick,
      domain=0.000000:0.500000
    ] {0.000000*(x-0.000000)+2.000000};
    \addplot[
      red,
      thick,
      domain=0.500000:0.666667
    ] {-2.000000*(x-0.500000)+2.000000};
    \addplot[
      red,
      thick,
      domain=0.666667:0.833333
    ] {-5.000000*(x-0.666667)+1.666667};
    \addplot[
      black,
      dashed,
      domain=0:0.833
    ] {x};
    \addplot[
      black,
      dashed
    ] coordinates {
      (0, 0)
      (0, 2.0)
    };
    \addplot[
      only marks,
      mark=*,
      mark size=3pt,
      red
    ] coordinates {
      (0.500000, 2.000000)
      (0.666667, 1.666667)
      (0.833333, 0.833333)
    };

    \addplot[
      only marks,
      mark=*,
      mark size=3pt,
      black
    ] coordinates {
      (0.000000, 0.000000)
      (0.000000, 2.000000) 
    };

    \node at (axis cs:0.375,1.1) [font=\huge] {many copies of $\bm{\mathcal{P}}$};
    
    \node at (axis cs:1.05,1.925) [font=\huge] {no copies of $\bm{\mathcal{P}}$};
    
    \node at (axis cs:0.150000,0.100000) [anchor=north east, xshift=-3pt, yshift=-3pt] {\Large{$O$}};
    \node at (axis cs:0.10000,1.975000) [anchor=south east, xshift=-3pt, yshift=3pt] {\Large{$Q_4$}};
    \node at (axis cs:0.500000,1.900000) [anchor=south west, xshift=3pt, yshift=3pt] {\Large{$Q_3$}};
    \node at (axis cs:0.666667,1.666667) [anchor=west, xshift=3pt] {\Large{$Q_2$}};
    \node at (axis cs:0.833333,0.833333) [anchor=west, xshift=3pt] {\Large{$Q_1$}};

    
  \end{axis}
\end{tikzpicture}

%% file: Revision.bbl
\begin{thebibliography}{10}

\bibitem{abraham2025probability}
R.~Abraham, J.-F. Delmas, and J.~Weibel.
\newblock Probability-graphons: {Limits} of large dense weighted graphs.
\newblock {\em Innovations in Graph Theory}, 2:25--117, 2025.

\bibitem{barak2019efficient}
B.~Barak, C.-N. Chou, Z.~Lei, T.~Schramm, and Y.~Sheng.
\newblock ({N}early) efficient algorithms for the graph matching problem on
  correlated random graphs.
\newblock {\em Advances in Neural Information Processing Systems}, 32, 2019.

\bibitem{barbour1989central}
A.~D. Barbour, M.~Karo{\'n}ski, and A.~Ruci{\'n}ski.
\newblock A central limit theorem for decomposable random variables with
  applications to random graphs.
\newblock {\em Journal of Combinatorial Theory, Series B}, 47(2):125--145,
  1989.

\bibitem{bentley2016multilayer}
B.~Bentley, R.~Branicky, C.~L. Barnes, Y.~L. Chew, E.~Yemini, E.~T. Bullmore,
  P.~E. V{\'e}rtes, and W.~R. Schafer.
\newblock The multilayer connectome of caenorhabditis elegans.
\newblock {\em {PL}o{S} computational Biology}, 12(12):e1005283, 2016.

\bibitem{bergermann2021multiplex}
J.~Bergermann and S.~Stoll.
\newblock Multiplex public transport networks: Centrality measures and
  vulnerability assessment.
\newblock {\em Applied Network Science}, 6(1):1--20, 2021.

\bibitem{bianconi2018multilayer}
G.~Bianconi.
\newblock {\em Multilayer networks: structure and function}.
\newblock Oxford university press, 2018.

\bibitem{bollobas1981threshold}
B.~Bollob{\'a}s.
\newblock Threshold functions for small subgraphs.
\newblock {\em Mathematical Proceedings of the Cambridge Philosophical
  Society}, 90(2):197--206, 1981.

\bibitem{bollobas1989subgraph}
B.~Bollob{\'a}s and J.~C. Wierman.
\newblock Subgraph counts and containment probabilities of balanced and
  unbalanced subgraphs in a large random graph.
\newblock {\em Annals of the New York Academy of Sciences}, 576(1):63--70,
  1989.

\bibitem{cardillo2013emergence}
A.~Cardillo, J.~G{\'o}mez-Gardenes, M.~Zanin, M.~Romance, D.~Papo, F.~d. Pozo,
  and S.~Boccaletti.
\newblock Emergence of network features from multiplexity.
\newblock {\em Scientific Reports}, 3(1):1344, 2013.

\bibitem{poissonapproximation}
S.~Chatterjee, P.~Diaconis, and E.~Meckes.
\newblock Exchangeable pairs and poisson approximation.
\newblock {\em Probability Surveys}, 2:64--106, 2005.

\bibitem{ChenGoldsteinShao2011}
L.~H.~Y. Chen, L.~Goldstein, and Q.-M. Shao.
\newblock {\em Normal Approximation by {S}tein's {M}ethod}.
\newblock Springer, Berlin, Heidelberg, 2011.

\bibitem{multiplexnetworkproperties}
E.~Cozzo, G.~F. De~Arruda, F.~A. Rodrigues, and Y.~Moreno.
\newblock {\em Multiplex networks: basic formalism and structural properties},
  volume~2.
\newblock Springer, 2018.

\bibitem{erdos1960evolution}
P.~Erd\H{o}s and A.~R{\'e}nyi.
\newblock On the evolution of random graphs.
\newblock {\em Publ. Math. Inst. Hung. Acad. Sci}, 5(1):17--60, 1960.

\bibitem{fienberg1985statistical}
S.~E. Fienberg, M.~M. Meyer, and S.~S. Wasserman.
\newblock Statistical analysis of multiple sociometric relations.
\newblock {\em Journal of the American Statistical Association},
  80(389):51--67, 1985.

\bibitem{multiplexnetworks}
A.~Ganguly and B.~B. Bhattacharya.
\newblock Multiplexons: Limits of multiplex networks.
\newblock {\em arXiv:2510.08639}, 2025.

\bibitem{huang2025correlation}
D.~Huang and P.~Yang.
\newblock Testing correlation in graphs by counting bounded degree motifs.
\newblock {\em arXiv:2510.25289}, 2025.

\bibitem{janson1990functional}
S.~Janson.
\newblock A functional limit theorem for random graphs with applications to
  subgraph count statistics.
\newblock {\em Random Structures \& Algorithms}, 1(1):15--37, 1990.

\bibitem{janson2011random}
S.~Janson, T.~{\L}uczak, and A.~Ruci\'{n}ski.
\newblock {\em Random graphs}, volume~45.
\newblock John Wiley \& Sons, 2000.

\bibitem{kahn2007thresholds}
J.~Kahn and G.~Kalai.
\newblock Thresholds and expectation thresholds.
\newblock {\em Combinatorics, Probability and Computing}, 16(3):495--502, 2007.

\bibitem{networksbiology}
N.~A. Kiani, D.~Gomez-Cabrero, and G.~Bianconi.
\newblock {\em Networks of Networks in Biology: Concepts, Tools and
  Applications}.
\newblock Cambridge University Press, 2021.

\bibitem{kirkpatrick2025shortest}
Y.~Kirkpatrick and V.~V. Williams.
\newblock Shortest paths in multimode graphs.
\newblock {\em arXiv:2506.22261}, 2025.

\bibitem{lyzinski2015graph}
V.~Lyzinski, D.~E. Fishkind, M.~Fiori, J.~T. Vogelstein, C.~E. Priebe, and
  G.~Sapiro.
\newblock Graph matching: Relax at your own risk.
\newblock {\em IEEE Transactions on Pattern Analysis and Machine Intelligence},
  38(1):60--73, 2016.

\bibitem{lyzinski2014seeded}
V.~Lyzinski, D.~E. Fishkind, and C.~E. Priebe.
\newblock Seeded graph matching for correlated {E}rd{\H{o}}s-{R}{\'e}nyi
  graphs.
\newblock {\em J. Mach. Learn. Res.}, 15(1):3513--3540, 2014.

\bibitem{mao2023random}
C.~Mao, Y.~Wu, J.~Xu, and S.~H. Yu.
\newblock Random graph matching at {O}tter's threshold via counting
  chandeliers.
\newblock In {\em Proceedings of the 55th Annual ACM Symposium on Theory of
  Computing}, pages 1345--1356, 2023.

\bibitem{mao2024testing}
C.~Mao, Y.~Wu, J.~Xu, and S.~H. Yu.
\newblock Testing network correlation efficiently via counting trees.
\newblock {\em The Annals of Statistics}, 52(6):2483--2505, 2024.

\bibitem{nowicki1989asymptotic}
K.~Nowicki.
\newblock Asymptotic normality of graph statistics.
\newblock {\em Journal of Statistical Planning and Inference}, 21(2):209--222,
  1989.

\bibitem{multiplexindex}
I.~M. Oliveira, L.~C. Carpi, and A.~P.~F. Atman.
\newblock The {M}ultiplex {E}fficiency {I}ndex: unveiling the {B}razilian air
  transportation multiplex network--{BATMN}.
\newblock {\em Scientific Reports}, 10(1):13339, 2020.

\bibitem{florentinegraph}
J.~F. Padgett and C.~K. Ansell.
\newblock Robust action and the rise of the {M}edici, 1400-1434.
\newblock {\em American Journal of Sociology}, 98(6):1259--1319, 1993.

\bibitem{conjectureproof}
J.~Park and H.~Pham.
\newblock A proof of the {K}ahn--{K}alai conjecture.
\newblock {\em Journal of the American Mathematical Society}, 37(1):235--243,
  2024.

\bibitem{pattison1999logit}
P.~Pattison and S.~Wasserman.
\newblock Logit models and logistic regressions for social networks: {II}.
  multivariate relations.
\newblock {\em British Journal of Mathematical and Statistical Psychology},
  52(2):169--193, 1999.

\bibitem{pedarsani2011privacy}
P.~Pedarsani and M.~Grossglauser.
\newblock On the privacy of anonymized networks.
\newblock In {\em 17th ACM SIGKDD International Conference on Knowledge
  Discovery and Data Mining}, pages 1235--1243, 2011.

\bibitem{perkins2025searching}
W.~Perkins.
\newblock Searching for (sharp) thresholds in random structures: Where are we
  now?
\newblock {\em Bulletin of the American Mathematical Society}, 62(1):113--143,
  2025.

\bibitem{rucinski1988small}
A.~Ruci{\'n}ski.
\newblock When are small subgraphs of a random graph normally distributed?
\newblock {\em Probability Theory and Related Fields}, 78(1):1--10, 1988.

\bibitem{rucinski1990survey}
A.~Ruci\'{n}ski.
\newblock Small subgraphs of random graphs--a survey.
\newblock In {\em Random graphs}, volume~87, pages 283--303, 1990.

\bibitem{strano2015multiplex}
E.~Strano, S.~Shai, S.~Dobson, and M.~Barthelemy.
\newblock Multiplex networks in metropolitan areas: generic features and local
  effects.
\newblock {\em Journal of The Royal Society Interface}, 12(111):20150651, 2015.

\bibitem{wasserman1994social}
S.~Wasserman and K.~Faust.
\newblock {\em Social Network Analysis: Methods and applications}.
\newblock Cambridge university press, 1994.

\bibitem{zucal2024probabilitygraphons}
G.~Zucal.
\newblock Probability graphons and {$P$}-variables: two equivalent viewpoints
  for dense weighted graph limits.
\newblock {\em arXiv:2408.07572}, 2024.

\end{thebibliography}
